\documentclass[reqno]{amsart}
\usepackage[matrix,arrow,curve,tips,frame]{xy}

\SelectTips{cm}{}

\usepackage{amsbsy,mathrsfs,latexsym,amssymb}

\theoremstyle{plain}
\newtheorem{theorem}{Theorem}[section]
\newtheorem{proposition}[theorem]{Proposition}
\newtheorem{corollary}[theorem]{Corollary}
\newtheorem{lemma}[theorem]{Lemma}

\theoremstyle{definition}
\newtheorem{definition}[theorem]{Definition}
\newtheorem{example}[theorem]{Example}

\newtheorem{assumption}[theorem]{Assumption}

\newtheorem{remark}[theorem]{Remark}

\newtheorem{notation}[theorem]{Notation}

\numberwithin{equation}{section}

\def\refeq#1{{\rm (\ref{#1})}}

\def\Colim#1#2{#1\star #2}

\def\kat#1{{\mathscr{#1}}}
\def\A{\kat{A}}
\def\B{\kat{B}}

\def\elts#1{{\mathsf{elts}}(#1)}
\def\Complex#1{{\mathsf{Complex}}(#1)}
\def\nComplex#1#2{{\mathsf{Complex}}_{#2}(#1)}

\def\opPresh#1{[#1^\op,\Set]}
\def\Presh#1{[#1,\Set]}
\def\Flat#1{{\mathsf{Flat}}(#1,\Set)}
\def\Lex#1{{\mathsf{Lex}}(#1,\Set)}

\def\zip{{\mathsf{zip}}}
\def\head{{\mathsf{head}}}
\def\tail{{\mathsf{tail}}}

\def\Set{{\mathsf{Set}}}
\def\Pos{{\mathsf{Pos}}}
\def\Inj{{\mathsf{Inj}}}
\def\Field{{\mathsf{Field}}}
\def\Lin{{\mathsf{Lin}}}

\def\Id{{\mathrm{Id}}}

\def\pr{{\mathrm{pr}}}
\def\beh{{\mathsf{beh}}}
\def\o{\cdot}
\def\op{{\mathit{op}}}
\def\fp{{\mathit{fp}}}

\def\phi{\varphi}
\def\tensor{\otimes}

\def\less{\sqsubseteq}
\def\strukt#1{\langle #1\rangle}
\def\after{\mathrel{;}}
\def\smash{\mathrel{\vee}}

\def\sol#1{#1^\dagger}

\def\komp#1#2{(#1_\bullet,#2_\bullet)}

\renewcommand{\to}{\longrightarrow}

\def\Lan#1#2{{\mathrm{Lan}}_{#1}{#2}}

\begin{document}
\title[Final Coalgebras in Accessible Categories]
      {Final Coalgebras in Accessible Categories}
\author{Panagis Karazeris}
\address{Department of Mathematics, University of Patras,
        Patras, Greece}
\email{pkarazer@math.upatras.gr}
\author{Apostolos Matzaris}
\address{Department of Mathematics, University of Patras,
        Patras, Greece}
\email{matzaris@master.math.upatras.gr}
\author{Ji\v{r}\'{\i} Velebil}
\address{Faculty of Electrical Engineering, Czech Technical University
         of Prague,
         Prague, Czech Republic}
\email{velebil@math.feld.cvut.cz}
\thanks{The authors acknowledge the support
        of the grant 
        MSM6840770014 of the Ministry of
        Education of 
        the Czech Republic that made their cooperation possible.
        Apostolos Matzaris also acknowledges the support of the 
        Carath\'{e}odory Basic Research Grant of the University of 
        Patras.}
\keywords{}
\subjclass{}
\date{29 May 2009}

\begin{abstract}
We give conditions on a finitary endofunctor of a finitely accessible 
category to admit a final coalgebra. Our conditions always apply
to the case of a finitary endofunctor of a locally finitely
presentable (l.f.p.) category and they bring an explicit construction 
of the final coalgebra in this case. On the other hand,
there are interesting examples of final coalgebras beyond the realm
of l.f.p. categories to which our results apply.
We rely on ideas developed by Tom Leinster for the study of 
self-similar objects in topology. 
\end{abstract}

\maketitle

\section{Introduction}
\label{sec:intro}

Coalgebras for an endofunctor (of, say, the category of sets)
are well-known to describe 
{\em systems of formal recursive equations\/}.
Such a system of equations then specifies a potentially
infinite ``computation'' and one is naturally interested 
in giving (uninterpreted) semantics to such a computation.
In fact, such semantics can be given by means of a coalgebra
again: this time by the {\em final\/} coalgebra for the given
endofunctor.

Let us give a simple example of that.

\begin{example}
\label{ex:simple}
Suppose that we fix a set $A$ and we want to consider
the set $A^\omega$ of infinite sequences of elements
of $A$, called {\em streams\/}. Moreover, we want
to define a function $\zip:A^\omega\times A^\omega\to A^\omega$
that ``zips up'' two streams, i.e., the equality
$$
\zip
\Bigl(
(a_0,a_1,a_2,\dots),
(b_0,b_1,b_2,\dots)
\Bigr)
=
(a_0,b_0,a_1,b_1,a_2,b_2,\dots)
$$
holds. 

One possible way of working with infinite expressions
like streams
is to introduce an additional approximation 
structure on the set of infinite expressions 
and to speak of an infinite expression     
as of a ``limit'' of its finite 
approximations, either in the sense of a complete partial order 
or of a complete metric
space, see~\cite{adj} and~\cite{america+rutten},
respectively. Such an approach may get rather
technical and the additional approximation structure
may seem rather arbitrary.

In fact, using the ideas of Calvin Elgot and his 
collaborators, see~\cite{elgot} and~\cite{ebt}, combined 
with a coalgebraic approach to systems of recursive
equations~\cite{rutten} and~\cite{aamv}, one may drop the additional
structure altogether and define solutions by 
{\em corecursion\/}, i.e., by means of a final coalgebra.

Clearly, the above zipping function can be specified
by a {\em system of recursive equations\/}
\begin{equation}
\label{eq:1}
\zip(a,b)
=
(\head(a),\zip(b,\tail(a)))
\end{equation}
one equation for each pair $a$, $b$ of streams, where
we have used the functions 
$\head(a_0,a_1,a_2,\dots)=a_0$ and 
$\tail(a_0,a_1,a_2,\dots)=(a_1,a_2,\dots)$.

In fact, the above system~\refeq{eq:1} of recursive equations 
can be encoded as a map
\begin{equation}
\label{eq:2}
e:A^\omega\times A^\omega \to A\times A^\omega \times A^\omega,
\quad
(a,b)\mapsto (\head(a),b,\tail(a))
\end{equation}
This means that we rewrote the system~\refeq{eq:1}
as a coalgebra and we will show now that a final coalgebra
gives its unique solution, namely the function $\zip$. 
To this end, we define first an endofunctor
$\Phi$ of the category of sets by the assignment 
$$
X\mapsto A\times X
$$
A {\em coalgebra\/} for $\Phi$
(with an {\em underlying set\/} $X$) 
is then any mapping $e:X\to\Phi X$, i.e., a mapping of the form
$$
e:X\to A\times X
$$

Suppose that a {\em final coalgebra\/} 
$$
\tau:TA\to A\times TA
$$
for $\Phi$ exists. Its finality means that for 
{\em any\/} coalgebra $c:Z\to A\times Z$
there exists a unique mapping $\sol{c}:Z\to TA$
such that the square
\begin{equation}
\label{eq:sol_square}
\vcenter{
\xymatrix{
Z
\ar[0,1]^-{c}
\ar[1,0]_{\sol{c}}
&
A\times Z
\ar[1,0]^{A\times\sol{c}}
\\
TA
\ar[0,1]_-{\tau}
&
A\times TA
}
}
\end{equation}
commutes. Moreover, it is well-known that the mapping
$\tau$ must be a bijection due to finality.
Luckily, in our case the final
coalgebra is well-known to exist and has the following
description: $TA$ is the set of all streams $A^\omega$
and the mapping $\tau$ sends $a\in A^\omega$ to the
pair $(\head(a),\tail(a))$.

If we instantiate the coalgebra $e$ 
from~\refeq{eq:2} for $c$ in the above square
and if we chase the elements of $A^\omega\times A^\omega$ 
around it, we see that the uniquely determined 
function $\sol{e}:A^\omega\times A^\omega\to A^\omega$
satisfies the recursive equation~\refeq{eq:1}. 
\end{example}

The reason for the existence of a final coalgebra for $\Phi$
is that both the category of sets and the endofunctor $\Phi$
are ``good enough'': the category of sets is locally finitely
presentable and the functor is finitary (we explain what that 
means in more detail below).

However, it is not the case that a final coalgebra exists
for every ``good enough'' functor: for example the identity
endofunctor of the category of sets and injections does not have
a final coalgebra for cardinality reasons.
Yet there are examples of interesting endofunctors of
``less good'' categories that still have a final coalgebra,
see, e.g., Example~\ref{ex:final_coalgebra_on_not_lfp_part1}
below. 

The important thing, however, is that our uniform description
of final coalgebras will be very reminiscent of streams:
the coalgebra structure of a final coalgebra is {\em always\/}
given by analogues of $\head$ and $\tail$ mappings from
the previous example.

\subsection*{The goals and organization of the paper}
In this paper we will focus on the existence of final coalgebras
for the class of finitary endofunctors of 
{\em finitely accessible categories\/}. Moreover, we will
give a concrete description of such coalgebras. From the
above it is clear how final coalgebras capture solutions
of recursive systems. 

We will make advantage of the fact that finitary endofunctors
of finitely accessible categories can be fully reconstructed
from {\em essentially small data\/}. In fact, finitary 
endofunctors can be replaced by
{\em flat modules\/} on the {\em small categories of finitely
presentable objects\/}. Such pairs
$$
\mbox{(small category,\ flat module)}
$$
will be called {\em self-similarity systems\/} and they 
fully encode the pattern of the recursive process in question.

We recall the concepts of finitary functors and finitely
accessible categories and the process of passing from
endofunctors to modules in Section~\ref{sec:prelim}.
 
In Section~\ref{sec:complexes} we
introduce the main tool of the paper --- the 
{\em category of complexes\/} for a (flat) module.
The category of complexes will then allow us to
give a concrete description of final coalgebras.

In Section~\ref{sec:SSC} we formulate a condition on the 
category of complexes that ensures that a final coalgebra
for the module in question exists, see Theorem~\ref{th:strong=>final}
below. As a byproduct we obtain, in Corollary~\ref{cor:barr}, 
a new proof of the well-known fact that every finitary
endofunctor of a {\em locally finitely presentable category\/}
has a final coalgebra. Moreover, we prove that the elements 
of the final coalgebra are essentially the complexes.  

Although the results of Section~\ref{sec:SSC} give a concrete
desription of the final coalgebra, the condition we give in
this section is rather strong. We devote Section~\ref{sec:WSC}
to a certain weakening of this condition. The weaker condition
on the category of complexes of the module yields a final
coalgebra as well but the module has to satisfy a certain
side condition of finiteness flavour. 

In some cases, one can prove that the conditions we give
are {\em necessary and sufficient\/} for the existence
of a final coalgebra. We devote Section~\ref{sec:final=>WSC}
to finding conditions on the endofunctor that ensure the
existence of such a characterization.

\subsection*{Related work}
This work is very much influenced by the work of Tom Leinster,
\cite{leinster1} and~\cite{leinster2} on self-similarity
in topology. In fact, Leinster works with categories
that are ``accessible'' for the notion of componentwise
filtered. 

Other descriptions of final coalgebras follow from 
the analysis of the final coalgebra sequence, see~\cite{adamek}.
However, this technique differs from ours.

\subsection*{Acknowledgements}  
We are grateful to Tom Leinster and Ji\v{r}\'{\i} Ad\'{a}mek 
for their valuable comments on earlier drafts of this paper.

\section{Preliminaries}
\label{sec:prelim}

In this preliminary section we introduce the notation and terminology
that we will use in the rest of the paper. Most of it is fairly
standard, we refer to books~\cite{ar} and~\cite{borceux}
for the material concerning finitely accessible categories
and finitary functors.

\subsection*{Coalgebras and final coalgebras}

We give a precise definition of (final) coalgebras,
see, e.g., \cite{rutten} for motivation and examples
of various coalgebras in the category of sets.

\begin{definition}
Suppose $\Phi:\kat{K}\to\kat{K}$ is any functor.
\begin{enumerate}
\item A {\em coalgebra for $\Phi$\/} is a morphism
      $e:X\to \Phi (X)$.
\item A {\em homomorphism of coalgebras\/} from
      $e:X\to\Phi(X)$ to $e':X'\to\Phi(X')$
      is a morphism $h:X\to X'$ making the following
      square 
      $$
      \xymatrix{
      X
      \ar[0,1]^-{e}
      \ar[1,0]_{h}
      &
      \Phi(X)
      \ar[1,0]^{\Phi(h)}
      \\
      X'
      \ar[0,1]_-{e'}
      &
      \Phi(X')
      }
      $$ 
      commutative.
\item A coalgebra $\tau:T\to\Phi(T)$ is called {\em final\/},
      if it is a terminal object of the category of coalgebras,
      i.e., if for every coalgebra $e:X\to\Phi(X)$ there is
      a unique morphism $\sol{e}:X\to T$ such that the square
      $$
      \xymatrix{
      X
      \ar[0,1]^-{e}
      \ar[1,0]_{\sol{e}}
      &
      \Phi(X)
      \ar[1,0]^{\Phi(\sol{e})}
      \\
      T
      \ar[0,1]_-{\tau}
      &
      \Phi(T)
      }
      $$ 
      commutes.
\end{enumerate}
\end{definition} 

\subsection*{Finitely accessible and locally finitely presentable categories}

Finitely accessible and locally finitely presentable categories are
those where every object can be reconstructed knowing its 
``finite parts''. This property has, for example, the category 
$\Set$ of sets and mappings, where a set $P$ is recognized as finite 
exactly when its hom-functor $\Set(P,{-}):\Set\to\Set$ preserves
colimits of a certain class --- the so-called {\em filtered\/}
colimits. 

A colimit of a general diagram $D:\kat{D}\to\kat{K}$ is called
{\em filtered\/}, provided that its scheme-category $\kat{D}$
is filtered. A category $\kat{D}$ is called {\em filtered\/} provided
that every finite subcategory of $\kat{D}$ admits a cocone.
In more elementary terms, filteredness of $\kat{D}$ can be expressed
equivalently by the following three properties:
\begin{enumerate}
\item The category $\kat{D}$ is nonempty.
\item Each pair $d_1$, $d_2$ of objects of $\kat{D}$ has an 
      ``upper bound'', i.e., there exists a cocone
      $$
      \xymatrixrowsep{.05pc}
      \xymatrix{
      d_1
      \ar @{.>} [1,1]
      &
      \\
      &
      d
      \\
      d_2
      \ar @{.>} [-1,1]
      & 
      }
      $$
      in $\kat{D}$.
\item Each parallel pair of morphisms in $\kat{D}$ can be
      ``coequalized'', i.e., for each parallel pair
      $$
      \xymatrix{
      d_1
      \ar @<.7ex> [0,1]
      \ar @<-.7ex> [0,1]
      &
      d_2
      }
      $$
      of morphisms in $\kat{D}$ there is a completion 
      to a commutative diagram of the form
      $$
      \xymatrix{
      d_1
      \ar @<.7ex> [0,1]
      \ar @<-.7ex> [0,1]
      &
      d_2
      \ar @{.>} [0,1]
      &
      d
      }
      $$      
      in $\kat{D}$.
\end{enumerate}
A category is $\kat{D}$ called {\em cofiltered\/} provided that 
the dual category $\kat{D}^\op$ is filtered.

An object $P$ of a category $\kat{K}$ is called {\em finitely presentable\/}
if the hom-functor $\kat{K}(P,{-}):\kat{K}\to\Set$ preserves
filtered colimits. 

\begin{definition}
\label{def:lfp}
A category $\kat{K}$ is called {\em finitely accessible\/} if it has
filtered colimits and if it contains a small subcategory consisting
of finitely presentable objects such that every object of $\kat{K}$
is a filtered colimit of these finitely presentable objects.

A cocomplete finitely accessible category is called 
{\em locally finitely presentable\/}.
\end{definition}

\begin{remark}
Locally finitely presentable categories were introduced by Peter
Gabriel and Friedrich Ulmer~\cite{gu}, finitely accessible
categories were introduced by Christian Lair~\cite{lair} under the name
{\em sketchable\/} categories. Tight connections of these
concepts to (infinitary) logic can be found in the book~\cite{mp},
the book~\cite{ar} deals with the connection of these concepts
to categories of structures. 
\end{remark}

\begin{example}
\label{ex:lfp}
\hfill
\begin{enumerate}
\item\label{item:sets}
      The category
      $
      \Set
      $
      of sets and mappings is locally finitely presentable. 
      The finitely presentable objects are exactly the finite sets.
\item Every variety of finitary algebras is a locally finitely
      presentable category. The finitely presentable objects
      are exactly the algebras that are presented by finitely
      many generators and finitely many equations in the sense
      of universal algebra. 
\item\label{item:injections}
      The category
      $
      \Inj
      $
      having sets as objects and injective maps
      as morphisms is a finitely accessible category
      that is not locally finitely presentable. The finitely presentable
      objects are exactly the finite sets.
\item\label{item:fields}
      Denote by
      $
      \Field
      $
      the category of fields and field homomorphisms. Then $\Field$
      is a finitely accessible category that is not locally finitely 
      presentable.
\item\label{item:linear}
      The category 
      $
      \Lin
      $
      of linear orders and monotone maps is finitely
      accessible but not locally finitely presentable. 
      The finitely presentable objects are exactly 
      the finite ordinals. 
\item\label{item:pos_01}
      Let $\Pos_{0,1}$ denote the following category:
      \begin{enumerate}
      \item Objects are posets having distinct top and bottom elements.
      \item Morphisms are monotone maps preserving top and bottom
            elements.
      \end{enumerate}
      Then $\Pos_{0,1}$ is a {\em Scott complete\/} category 
      in the sense of Ji\v{r}\'{\i} Ad\'{a}mek~\cite{adamek_scc}: 
      it is finitely accessible and every small diagram in $\Pos_{0,1}$
      that has a cocone, has a colimit. 
 
      Scott complete categories are therefore
      ``not far away'' from being cocomplete and thus
      locally finitely presentable.

      However, $\Pos_{0,1}$ is not locally finitely presentable
      since it lacks a terminal object. Finitely presentable
      objects in $\Pos_{0,1}$ are exactly the finite posets having
      distinct bottom and top elements.
\item The category of topological spaces and continuous maps
      is not finitely accessible. Although this category has
      filtered (in fact, all) colimits, the only finitely
      presentable objects are finite discrete topological spaces 
      and these do not suffice for reconstruction of 
      a general topological space.
\end{enumerate}
Of course, more examples of ``everyday-life'' finitely accessible
categories can be found in the literature, see, e.g., 
papers~\cite{diers1} and~\cite{diers2} by Yves Diers.
\end{example}

Every finitely accessible category $\kat{K}$ is equivalent
to a category of the form
$$
\Flat{\A}
$$
(where $\A$ is a small category) 
that consists of all {\em flat functors\/} $X:\A\to\Set$ 
and all natural transformations between them. 

A functor $X:\A\to\Set$ is called {\em flat\/} if its 
{\em category of elements\/} $\elts{X}$ is cofiltered. 
The category $\elts{X}$
has pairs $(x,a)$ with $x\in Xa$ as objects and as morphisms
from $(x,a)$ to $(x',a')$ those morphisms $f:a\to a'$ in $\A$ 
with the property that $Xf(x)=x'$.

Flat functors
$X$ can be characterized by any of the following equivalent
conditions:
\begin{enumerate}
\item The functor $X:\A\to\Set$ is a filtered colimit of 
      representable functors.
\item The left Kan extension $\Lan{Y}{X}:\opPresh{\A}\to\Set$
      of $X:\A\to\Set$ along the Yoneda embedding
      $Y:\A\to\opPresh{\A}$
      preserves finite limits. 
\end{enumerate}

In case when $\kat{K}$ is locally finitely presentable one can 
prove that $\kat{K}$ is equivalent to the category
$$
\Lex{\A}
$$
of all finite-limits-preserving functors on a small finitely
complete category $\A$. In fact, the flat functors are exactly
the finite-limits-preserving ones in this case.

\begin{example}
\label{ex:set=flat}
In this example we show how to express $\Set$ as a category 
of flat functors.
Denote by $E:\Set_\fp\to\Set$ the full dense inclusion of an 
essentially small category of finite sets. In fact, in this example, 
we choose as a representative set of finitely presentable objects the
set of finite ordinals. 

The correspondence
$$
X\mapsto \Set(E{-},X)
$$ 
then provides us with an equivalence
$$
\Set
\simeq
\Flat{\Set_\fp^\op}
=
\Lex{\Set_\fp^\op}
$$
of categories. The slogan behind this correspondence is
the following one:
\begin{itemize}
\item[] Instead of describing a set $X$ by means of its
        elements $x\in X$ (as we do in $\Set$), we describe
        a set by ``generalized elements'' of the form 
        $n\to X$, where $n$ is a finite ordinal.
\end{itemize}
Thus, a set $X$ now ``varies in time'': the hom-set 
$\Set(n,X)$ is the ``value'' of $X$ at ``time'' $n$.
\end{example}

\begin{remark}
\label{rem:fa=flat}
The above example is an instance of a general fact:
every finitely accessible category $\kat{K}$
is equivalent to $\Flat{\kat{K}_\fp^\op}$, where 
$E:\kat{K}_\fp\to\kat{K}$
denotes the full inclusion of the essentially small 
subcategory consisting of finitely presentable objects. 

The equivalence works as follows: the flat functor 
$X:\kat{K}_\fp^\op\to\Set$ is sent to the object
$$
\Colim{X}{E}
$$
which is a colimit of $E$ weighted by $X$. Such a colimit
is defined as an object $\Colim{X}{E}$ together with an
isomorphism
$$
\kat{K}(\Colim{X}{E},Z)
\cong
[\kat{K}_\fp^\op,\Set](X,\kat{K}(E{-},Z))
$$
natural in $Z$. The above colimit can be considered to be
an ``ordinary'' colimit of the diagram of elements
of $X$:
$$
x\in Xa \mapsto Ea
$$
This explains the weight terminology: every $Ea$ is going to be counted
``$Xa$-many times'' in the colimit $\Colim{X}{E}$. 
See~\cite{borceux} for more details.
\end{remark}

\subsection*{Flat modules}

On finitely accessible categories there is class of functors
that can be fully reconstructed by knowing their values
on ``finite parts''. An example is the finite-powerset
endofunctor
$$
P_{\mathit{fin}}:
X
\mapsto
\{ S\mid S\subseteq X, \mbox{ $S$ is finite } \}
$$
of the category of sets. Such endofunctors can be characterized
as exactly those {\em preserving\/} filtered colimits.

\begin{definition}
\label{def:finitary}
A functor $\Phi:\kat{K}\to\kat{L}$ between finitely accessible
categories is called {\em finitary\/} if it preserves filtered
colimits.
\end{definition} 

By the above considerations, every finitary
endofunctor $\Phi:\kat{K}\to\kat{K}$ of a finitely 
accessible category $\kat{K}$ can be considered, 
to within equivalence, as a finitary endofunctor 
$$
\Phi:\Flat{\A}\to\Flat{\A}
$$
Since the full embedding $\A^\op\to\Flat{\A}$ exhibits
$\Flat{\A}$ as a free cocompletion of $\A^\op$ w.r.t. filtered
colimits, we can then reconstruct $\Phi$ from a mere functor
$$
M_\Phi:\A^\op\to\Flat{\A}
$$
(no preservation properties) by means of filtered colimits.

The latter functor can be identified with a functor of the
form $M_\Phi:\A^\op\times\A\to\Set$ with the property that
every $M_\Phi(a,{-}):\A\to\Set$ is flat. Such functors of
two variables (without the extra flatness property) are commonly
called {\em modules\/}. We will give the extra property
a name.

\begin{definition}
\label{def:flat_module}
A {\em module\/} 
$
M:\A
\xymatrixcolsep{1.5pc}\xymatrix{\ar[0,1]|-{\object @{/}}&}
\B
$
from a small category $\A$ to a small category $\B$ is a functor
$M:\A^\op\times\B\to\Set$. Given two such modules, $M$ and $N$,
a {\em module morphism\/}
$M\to N$ is a natural transformation between the respective 
functors.

A module $M$ as above is called {\em flat\/}
if every partial functor $M(a,{-}):\B\to\Set$ is a flat
functor in the usual sense. 
\end{definition}

\begin{remark}
\label{rem:action}
The above module terminology makes perfect sense if we denote
an element $m\in M(a,b)$ by an arrow
$$
\xymatrix@1{
a
\ar[0,1]|-{\object @{/}}^-{m}
&
b
}
$$
and think of it as of a ``vector'' on which the categories
$\A$ and $\B$ can act by means of their morphisms (``scalars''):
\begin{enumerate}
\item Given $f:a'\to a$ in $\A$, then 
      $$
      \xymatrix@1{
      a'
      \ar[0,1]^-{f}
      &
      a
      \ar[0,1]|-{\object @{/}}^-{m}
      &
      b
      }
      $$ 
      denotes the element $M(f,b)(m)\in M(a',b)$.

      Had we denoted such an action by $m@f$, then
      it is obvious that equations 
      $m@(f\o f')=(m@f)@f'$ and $m@ 1_{a}=m$ hold ---
      something that we know from classical module
      theory. 
\item Given $g:b\to b'$ in $\B$, then 
      $$
      \xymatrix@1{
      a
      \ar[0,1]|-{\object @{/}}^-{m}
      &
      b
      \ar[0,1]^-{g}
      &
      b'
      }
      $$ 
      denotes the element $M(a,g)(m)\in M(a,b')$.
\item Functoriality of $M$ gives an unambiguous meaning
      to diagrams of the form
      $$
      \xymatrix@1{
      a'
      \ar[0,1]^-{f}
      &
      a
      \ar[0,1]|-{\object @{/}}^-{m}
      &
      b
      \ar[0,1]^-{g}
      &
      b'
      }
      $$
\item We also extend the notion of commutative diagrams. For example,
      by saying that the following square
      $$
      \xymatrix{
      a
      \ar[0,1]|-{\object @{/}}^-{m}
      \ar[1,0]_{f}
      &
      b
      \ar[1,0]^{g}
      \\
      a'
      \ar[0,1]|-{\object @{/}}_-{m'}
      &
      b'
      }
      $$
      commutes we mean that the equality $m'@ f=g@m$ holds.
\end{enumerate}
\end{remark}

\begin{remark}
\label{rem:flatness_of_M}
The broken arrow notation also allows us to formulate flatness
of a module
$
M:\A
\xymatrixcolsep{1.5pc}\xymatrix{\ar[0,1]|-{\object @{/}}&}
\B
$
in elementary terms. Namely, for every $a$ in $\A$ the 
following three conditions must be satisfied:
\begin{enumerate}
\item There is a broken arrow
      $$
      \xymatrix@1{
      a
      \ar[0,1]|-{\object @{/}}^-{m}
      &
      b
      }
      $$ 
      for some $b$ in $\B$.
\item For any two broken arrows
      $$
      \xymatrixrowsep{.25pc}
      \xymatrix{
      &
      b_1
      \\
      a
      \ar[-1,1]|-{\object @{/}}^-{m_1}
      \ar[1,1]|-{\object @{/}}_-{m_2}
      &
      \\
      &
      b_2
      }
      $$       
      there is a commutative diagram
      $$
      \xymatrixrowsep{1.5pc}
      \xymatrix{
      &
      &
      b_1
      \\
      a
      \ar[-1,2]|-{\object @{/}}^-{m_1}
      \ar @{.>} [0,1]|-{\object @{/}}^-{m}
      \ar[1,2]|-{\object @{/}}_-{m_2}
      &
      b
      \ar @{.>} [-1,1]_{f_1}
      \ar @{.>} [1,1]^{f_2}
      &
      \\
      &
      &
      b_2
      }
      $$
\item For every commutative diagram
      $$
      \xymatrixrowsep{.25pc}
      \xymatrix{
      &
      b_1
      \ar @<-.5ex> [2,0]_{u}
      \ar @<.5ex> [2,0]^{v}
      \\
      a
      \ar[-1,1]|-{\object @{/}}^-{m_1}
      \ar[1,1]|-{\object @{/}}_-{m_2}
      &
      \\
      &
      b_2
      }
      $$    
      there is a commutative diagram
      $$
      \xymatrixrowsep{1.5pc}
      \xymatrix{
      &
      b
      \ar @{.>} [1,0]^{f}
      \\
      a
      \ar @{.>} [-1,1]|-{\object @{/}}^-{m}
      \ar[0,1]|-{\object @{/}}^-{m_1}
      \ar[1,1]|-{\object @{/}}_-{m_2}
      &
      b_1
      \ar @<-.5ex> [1,0]_{u}
      \ar @<.5ex> [1,0]^{v}
      \\
      &
      b_2
      }
      $$    
\end{enumerate}
\end{remark}

\begin{example}
\label{ex:phi_on_sets}
In this example we show how the finitary endofunctor
$$
X\mapsto X\times X+A
$$
of the locally finitely presentable category $\Set$ can be viewed
as a flat module. 

In this sense, we identify the endofunctor 
$X\mapsto X\times X+A$ of $\Set$ with the endofunctor
$$
\Phi:\Set(E{-},X)\mapsto\Set(E{-},X\times X)+\Set(E{-},A)
$$
of $\Flat{\Set_\fp^\op}$. The corresponding flat
module
$$
M:
\Set_\fp^\op
\xymatrixcolsep{1.5pc}\xymatrix{\ar[0,1]|-{\object @{/}}&}
\Set_\fp^\op
$$
then has values
$$
M(a,b)=\Set_\fp (b,a\times a)+\Set (b,A)
$$
at finite ordinals $a$, $b$. 
\end{example}

The above resemblance to classical module 
theory\footnote{The resemblance can be made precise by passing to
enriched category theory, see~\cite{borceux}.} 
can be pushed further: modules can 
composed by ``tensoring'' them.

\begin{definition}
Suppose
$
M:\A
\xymatrixcolsep{1.5pc}\xymatrix{\ar[0,1]|-{\object @{/}}&}
\B
$
and
$
N:\B
\xymatrixcolsep{1.5pc}\xymatrix{\ar[0,1]|-{\object @{/}}&}
\kat{C}
$
are modules. By
$$
N\tensor M:\A
\xymatrixcolsep{1.5pc}\xymatrix{\ar[0,1]|-{\object @{/}}&}
\kat{C}
$$
we denote their {\em composition\/} which is defined objectwise
by means of a coend
$$
\Bigl(N\tensor M\Bigr)(a,c)
=
\int^b N(b,c)\times M(a,b)
$$
\end{definition}

\begin{remark}
\label{rem:coends}
A coend is a special kind of colimit. The elements
of $\Bigl(N\tensor M\Bigr)(a,c)$ are equivalence classes.
A typical element of $\Bigl(N\tensor M\Bigr)(a,c)$ is
an equivalence class $[(n,m)]$ represented by a pair
$(n,m)\in N(b,c)\times M(a,b)$ where the equivalence is 
generated by requiring the pairs
$$
(n,f@m)
\quad
\mbox{and}
\quad
(n@f,m)
$$
to be equivalent, where $n$, $f$ and $m$ are as follows:
$$
\xymatrix@1{
a
\ar[0,1]|-{\object @{/}}^-{m}
&
b
\ar[0,1]^-{f}
&
b'
\ar[0,1]|-{\object @{/}}^-{n}
&
c
}
$$
Above, we denoted the actions of $M$ and $N$ by the same
symbols, not to make the notation heavy.
\end{remark}

It is well-known (see~\cite{borceux}) that the above composition 
organizes modules into a {\em bicategory\/}: 
the composition is associative only up
to a coherent isomorphism and the {\em identity module\/}
$
\A:\A
\xymatrixcolsep{1.5pc}\xymatrix{\ar[0,1]|-{\object @{/}}&}
\A
$,
sending $(a',a)$ to the hom-set $\A(a',a)$, serves as a unit
only up to a coherent isomorphism. The following result is then
easy to prove.

\begin{lemma}
Every identity module is flat and composition of flat modules
is a flat module.
\end{lemma}

\begin{remark}
The above composition of modules makes one to attempt
to draw diagrams such as
$$
\xymatrix@1{
a_2
\ar[0,1]|-{\object @{/}}^-{m_2}
&
a_1
\ar[0,1]|-{\object @{/}}^-{m_1}
&
a_0
}
$$
for elements $m_1\in M(a_1,a_0)$, $m_2\in M(a_2,a_1)$
of a module 
$
M:\A
\xymatrixcolsep{1.5pc}\xymatrix{\ar[0,1]|-{\object @{/}}&}
\A
$.
Such diagrams are, however, to be 
considered only formally --- {\em we never compose
two ``broken'' arrows\/}.
\end{remark}

The tensor notation from the above paragraphs allows us to pass
from endofunctors to modules completely. 

Observe that any flat functor $X:\A\to\Set$ can be considered
as a flat module
$
X:{\bf 1}
\xymatrixcolsep{1.5pc}\xymatrix{\ar[0,1]|-{\object @{/}}&}
\A
$
where ${\bf 1}$ denotes the one-morphism category.

Then, given a flat module
$
M:\A
\xymatrixcolsep{1.5pc}\xymatrix{\ar[0,1]|-{\object @{/}}&}
\A
$,
the assignment $X\mapsto M\tensor X$ defines a finitary
endofunctor of $\Flat{\A}$. 

In fact, every finitary endofunctor $\Phi$ of $\Flat{\A}$ arises 
in the above way: construct the flat module $M_\Phi$ as above, then
there is an isomorphism 
$$
\Phi
\cong
M_\Phi\tensor {-}
$$ 
of functors.

\section{The category of complexes and self-similarity systems}
\label{sec:complexes}
Formal chains of ``broken arrows'' will be the 
main tool of the rest of the paper. We define a category
of such chains (this definition comes from the 
paper~\cite{leinster1} of Tom Leinster).

\begin{assumption}
In the rest of the paper,
$$
M:\A
\xymatrixcolsep{1.5pc}\xymatrix{\ar[0,1]|-{\object @{/}}&}
\A
$$
denotes a flat module on a small category $\A$.
The pair $(\A,M)$ is called a
{\em self-similarity system\/}. 
\end{assumption}

\begin{remark}
The terminology {\em self-similarity system\/} is due
to Tom Leinster~\cite{leinster1} and has its origin
in the intention to study (topological) spaces that 
are self-similar. Since we refer to~\cite{leinster1}
below, we keep the terminology, although our motivation
is different.
\end{remark}

\begin{definition}
\label{def:complex}
Given a (flat) module $M$, the category
$$
\Complex{M}
$$
of {\em $M$-complexes and their morphisms\/} is defined as follows:
\begin{enumerate}
\item Objects, called {\em $M$-complexes\/}, are countable chains
      of the form
      $$
      \xymatrix{
      \dots
      \ar[0,1]|-{\object @{/}}^-{m_3}
      &
      a_2
      \ar[0,1]|-{\object @{/}}^-{m_2}
      &
      a_1
      \ar[0,1]|-{\object @{/}}^-{m_1}
      &
      a_0
      }
      $$
      A single complex as above will be denoted by
      $\komp{a}{m}$ for short.
\item {\em Morphisms from $\komp{a}{m}$ to $\komp{a'}{m'}$\/}
      are sequences $f_n:a_n\to a'_n$, denoted by $(f_\bullet)$,
      such that all squares in the following diagram
      $$
      \xymatrix{
      \dots
      \ar[0,1]|-{\object @{/}}^-{m_3}
      &
      a_2
      \ar[0,1]|-{\object @{/}}^-{m_2}
      \ar[1,0]^{f_2}
      &
      a_1
      \ar[0,1]|-{\object @{/}}^-{m_1}
      \ar[1,0]^{f_1}
      &
      a_0
      \ar[1,0]^{f_0}
      \\
      \dots
      \ar[0,1]|-{\object @{/}}_-{m'_3}
      &
      a'_2
      \ar[0,1]|-{\object @{/}}_-{m'_2}
      &
      a'_1
      \ar[0,1]|-{\object @{/}}_-{m'_1}
      &
      a'_0
      }
      $$
      commute. 
\end{enumerate}
For $n\geq 0$, we denote by
$$
\nComplex{M}{n}
$$
the category of {\em $n$-truncated $M$-complexes\/}.
Its objects are finite chains
$$
\xymatrix{
a_n      
\ar[0,1]|-{\object @{/}}^-{m_n}
&
a_{n-1}
\ar[0,1]|-{\object @{/}}
&
\dots
\ar[0,1]|-{\object @{/}}
&
a_2
\ar[0,1]|-{\object @{/}}^-{m_2}
&
a_1
\ar[0,1]|-{\object @{/}}^-{m_1}
&
a_0
}
$$
and the morphisms of $\nComplex{M}{n}$ are defined in the
obvious way.

The obvious truncation functors are denoted by
$$
\pr_n:\Complex{M}\to\nComplex{M}{n},
\quad
n\geq 0
$$
Observe that $\nComplex{M}{0}=\A$.
\end{definition}

\begin{example}
\label{ex:complex}
Recall the flat module $M$ of Example~\ref{ex:phi_on_sets}
that corresponds to the finitary endofunctor $X\mapsto X\times X+A$
of sets. 
 
An $M$-complex 
$$
\xymatrix{
\dots
\ar[0,1]|-{\object @{/}}^-{m_3}
&
a_2
\ar[0,1]|-{\object @{/}}^-{m_2}
&
a_1
\ar[0,1]|-{\object @{/}}^-{m_1}
&
a_0
}
$$
can be identified with a ``binary tree'' of maps of the form
$$
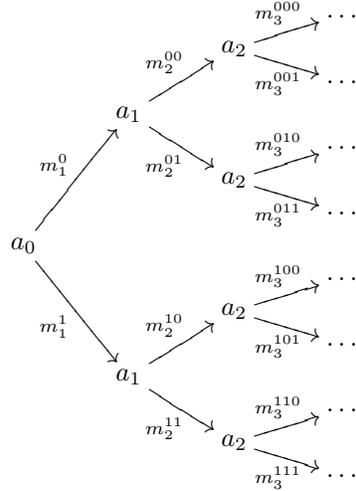
\xymatrixrowsep{.25pc}
\xymatrix{
&
&
&
\dots
\\
&
&
a_2
\ar[-1,1]^{m^{000}_{3}}
\ar[1,1]_{m^{001}_{3}}
&
\\
&
&
&
\dots
\\
&
a_1
\ar[-2,1]^{m^{00}_{2}}
\ar[2,1]_{m^{01}_{2}}
&
&
\\
&
&
&
\dots
\\
&
&
a_2
\ar[-1,1]^{m^{010}_{3}}
\ar[1,1]_{m^{011}_{3}}
&
\\
&
&
&
\dots
\\
a_0
\ar[-4,1]^{m^{0}_{1}}
\ar[4,1]_{m^{1}_{1}}
&
&
&
\\
&
&
&
\dots
\\
&
&
a_2
\ar[-1,1]^{m^{100}_{3}}
\ar[1,1]_{m^{101}_{3}}
&
\\
&
&
&
\dots
\\
&
a_1
\ar[-2,1]^{m^{10}_{2}}
\ar[2,1]_{m^{11}_{2}}
&
&
\\
&
&
&
\dots
\\
&
&
a_2
\ar[-1,1]^{m^{110}_{3}}
\ar[1,1]_{m^{111}_{3}}
&
\\
&
&
&
\dots
}
$$
where each path is either infinite or it ends with a generalized
element of $A$.
\end{example}

\begin{remark}
\label{rem:complexes=simple}
The description of complexes is particularly simple if one starts
with a finitely accessible category $\kat{K}$ and a finitary
endofunctor $\Phi:\kat{K}\to\kat{K}$. Then a complex (for the 
module corresponding to $\Phi$) is just a sequence
$$
\xymatrix{
a_0
\ar[0,1]^-{m_1}
&
\Phi(a_1)
,
&
a_1
\ar[0,1]^-{m_2}
&
\Phi(a_2)
,
&
a_2
\ar[0,1]^-{m_3}
&
\Phi(a_3)
,
&
\dots
}
$$
of morphisms in $\kat{K}$, where all the objects $a_0$, $a_1$,
$a_2$, \dots are finitely presentable.

And morphisms of complexes are just sequences of morphisms
in $\kat{K}$ making the obvious squares commutative:
$$
\xymatrix{
a_0
\ar[0,1]^-{m_1}
\ar[1,0]_{f_0}
&
\Phi(a_1)
\ar[1,0]^{\Phi (f_1)}
&
a_1
\ar[0,1]^-{m_2}
\ar[1,0]_{f_1}
&
\Phi(a_2)
\ar[1,0]^{\Phi (f_2)}
&
a_2
\ar[0,1]^-{m_3}
\ar[1,0]_{f_2}
&
\Phi(a_3)
\ar[1,0]^{\Phi (f_3)}
&
\ar@{} [1,1]|{\objectstyle\dots}
&
\\
a'_0
\ar[0,1]_-{m'_1}
&
\Phi(a'_1)
&
a'_1
\ar[0,1]_-{m'_2}
&
\Phi(a'_2)
&
a'_2
\ar[0,1]_-{m'_3}
&
\Phi(a'_3)
&
&
}
$$
In fact, a complex seen in this way is a ``finitary bit''
of a general coalgebra in the following sense:
start with a coalgebra $c:X\to \Phi(X)$
and a morphism $f_0:a_0\to X$, where $a_0$ is finitely
presentable. Due to finitarity of $\Phi$, 
the composite $c\o f_0:a_0\to\Phi(X)$ factors
through $\Phi(f_1):\Phi(a_1)\to\Phi(X)$ where $a_1$
is finitely presentable. The factorizing map
$m_1:a_0\to\Phi(a_1)$ is then the germ of a complex:
proceed with $f_1:a_1\to X$ to obtain $m_2:a_1\to\Phi(a_2)$,
etc. 
\end{remark}

\section{The Strong Solvability Condition}
\label{sec:SSC}

The Strong Solvability Condition on a self-similarity system $(\A,M)$ 
will give us a final coalgebra for the finitary functor
$$
M\tensor {-}:\Flat{\A}\to\Flat{\A}
$$
almost ``for free''. The condition asserts that there
is a certain filtered diagram of representables in 
$\Flat{\A}$.
The carrier of the final coalgebra is simply its
colimit, see Theorem~\ref{th:strong=>final} below. 

Although the condition is rather strong and
hard to verify in practice (and we will seek a weaker 
one in next section), it is trivially satisfied in the 
realm of locally finitely presentable categories. Hence
the technique of the current section enables us to give 
a uniform description of final
coalgebras for finitary endofunctors of locally finitely
presentable categories, see Corollary~\ref{cor:barr}.

Most of the results of this section are reformulations
of things proved in~\cite{leinster1} by Tom Leinster into
our setting. 

We give first an example of a finitely accessible category $\kat{K}$
that is not locally finitely presentable
and a finitary endofunctor $\Phi:\kat{K}\to\kat{K}$ that admits
a final coalgebra.

\begin{example}
\label{ex:final_coalgebra_on_not_lfp_part1}
Recall the category $\Pos_{0,1}$ of all posets having distinct top
and bottom and all monotone maps preserving top and bottom
of Example~\ref{ex:lfp}\refeq{item:pos_01}. Recall that $\Pos_{0,1}$
is finitely accessible but not locally finitely presentable.

It has been shown by Peter Freyd~\cite{freyd} that there is 
a finitary endofunctor $\Phi$ of $\Pos_{0,1}$ whose final
coalgebra gives the unit interval $[0,1]$.

The functor $\Phi:\Pos_{0,1}\to\Pos_{0,1}$ sends $X$ to
the {\em smash coproduct}
$$
X\smash X
$$
of $X$ with itself that is defined as follows: put one
copy of $X$ on top of the other one and glue the copies 
together by identifying top and bottom. More formally,
$X\smash X$ is the subposet of $X\times X$ consisting
of pairs $(x,0)$ or $(1,y)$. The pairs $(x,0)$
are going to be called living in the {\em left-hand\/} copy
of $X$ and the pairs of the form $(1,y)$ as living
in the {\em right-hand\/} copy.

Clearly, given a coalgebra $e:X\to X\smash X$ and $x\in X$,
one can produce at least one infinite sequence 
$$
x_1 x_2 x_3 \dots
$$
of 0's and 1's as follows: look at $e(x)$ and put $x_1=0$
if it is in the left-hand copy of $X$, put $x_1=1$ otherwise.
Then regard $e(x)$ as an element of $X$ again, apply
$e$ to it to produce $x_2$, etc.

One needs to show that the binary expansion
$\sol{e}(x)=0.x_1 x_2 x_3\dots$ so obtained can be used 
to define a map $\sol{e}:X\to [0,1]$ in a clash-free way
(i.e., regardless of the fact that sometimes we may have
a choice in defining $x_k=0$ or $x_k=1$). Moreover,
the above map $\sol{e}$ is then a witness that the
coalgebra
$$
t:[0,1]\to [0,1]\smash [0,1]
$$
where $[0,1]$ denotes the closed unit interval
with the usual order and $t$ given by putting
$t(x)=(2x,0)$ for $0\leq x\leq 1/2$ and 
$t(x)=(1,2x-1)$ otherwise, is a final coalgebra
for $\Phi$.

See~\cite{freyd} for more details on finality
of $[0,1]$ and see Example~\ref{ex:final_coalgebra_on_not_lfp_part3}
below that the description of a final
coalgebra for $\Phi$ given by our theory 
will provide us with the unit interval canonically.
\end{example}

We introduce now a condition on a self-similarity system 
$(\A,M)$ that will ensure the existence of a final coalgebra.

\begin{definition}
\label{def:SSC}
We say that $(\A,M)$ satisfies the
{\em Strong Solvability Condition\/} if the category
$\Complex{M}$ is cofiltered.
\end{definition}

\begin{remark}
\label{rem:SSC=>flat}
The Strong Solvability Condition implies that the diagram
$$
\xymatrix{
\Bigl(\Complex{M}\Bigr)^\op
\ar[0,1]^-{\pr_0^\op}
&
\A^\op
\ar[0,1]^-{Y}
&
\Presh{\A}
}
$$
of representables is filtered. Its colimit 
(a flat functor!) is going to be the carrier of the final
coalgebra for $M\tensor {-}$, see Theorem~\ref{th:strong=>final}
below.
\end{remark}

\begin{remark}
\label{rem:SSC=elementary}
In elementary terms, the Strong Solvability Condition
says that the following three conditions hold:
\begin{enumerate}
\item The category $\Complex{M}$ is nonempty.
\item For every pair $\komp{a}{m}$,
      $\komp{a'}{m'}$ in $\Complex{M}$ there is a span
      $$
      \xymatrixrowsep{.5pc}
      \xymatrix{
      &
      \komp{a}{m}
      \\
      \komp{b}{n}
      \ar[-1,1]^{(f_\bullet)}
      \ar[1,1]_{(f'_\bullet)}
      &
      \\
      &
      \komp{a'}{m'}
      }
      $$
      in $\Complex{M}$.
\item For every parallel pair of the form
      $$
      \xymatrix{
      \komp{a}{m}
      \ar @<1ex> [0,1]^-{(u_\bullet)}
      \ar @<-1ex> [0,1]_-{(v_\bullet)}
      &
      \komp{a'}{m'}
      }
      $$
      in $\Complex{M}$
      there is a fork
      $$
      \xymatrix{
      \komp{b}{n}
      \ar[0,1]^-{(f_\bullet)}
      &
      \komp{a}{m}
      \ar @<1ex> [0,1]^-{(u_\bullet)}
      \ar @<-1ex> [0,1]_-{(v_\bullet)}
      &
      \komp{a'}{m'}
      }
      $$
      in $\Complex{M}$.
\end{enumerate}
\end{remark}

\begin{example}
\label{ex:final_coalgebra_on_not_lfp_part2}
(Continuation of Example~\ref{ex:final_coalgebra_on_not_lfp_part1}.)

We show that the self-similarity system $(\A,M)$ corresponding
to the functor $\Phi:\Pos_{0,1}\to\Pos_{0,1}$ of 
Example~\ref{ex:final_coalgebra_on_not_lfp_part1}
satisfies the Strong Solvability Condition.

Recall that $M$ is defined
as
$$
M(a,b)
=
\Pos_{0,1}(b,a\smash a)
$$
where the posets $a$, $b$ are finite
(having distinct bottom and top).

A complex $\komp{a}{m}$ is therefore
a chain
$$
m_1:a_0\to a_1\smash a_1,
\quad
m_2:a_1\to a_2\smash a_2,
\quad
\dots,
\quad
m_i:a_i\to a_{i+1}\smash a_{i+1},
\quad
\dots
$$
of morphisms in $\Pos_{0,1}$.

We have to show that $\Complex{M}$ is cofiltered
and we will use the elementary description of
complexes of Remark~\ref{rem:complexes=simple} and
the elementary description of cofilteredness of 
Remark~\ref{rem:SSC=elementary}:
\begin{enumerate}
\item $\Complex{M}$ is nonempty.

      Let $a_i=2$, the two-element chain, for every $i\geq 0$
      and, for all $i\geq 0$, let $m_i:a_i\to a_{i+1}\smash a_{i+1}$ 
      be the unique morphism in $\Pos_{0,1}$. This defines a complex.
\item $\Complex{M}$ has cones for two-element discrete diagrams.

      Suppose $\komp{a}{m}$ and $\komp{a'}{m'}$ are given.
      Hence we have chains
      $$
      m_1:a_0\to a_1\smash a_1,
      \quad
      m_2:a_1\to a_2\smash a_2,
      \quad
      \dots,
      \quad
      m_i:a_i\to a_{i+1}\smash a_{i+1},
      \quad
      \dots
      $$
      and
      $$
      m'_1:a'_0\to a'_1\smash a'_1,
      \quad
      m'_2:a'_1\to a'_2\smash a'_2,
      \quad
      \dots,
      \quad
      m'_i:a'_i\to a'_{i+1}\smash a'_{i+1},
      \quad
      \dots
      $$

      Since every pair $a_i$, $a'_i$ has a cocone in
      $(\Pos_{0,1})_\fp$, every pair $a_i$, $a'_i$ 
      has a coproduct $a_i+a'_i$ in $(\Pos_{0,1})_\fp$ due to
      Scott-completeness of $\Pos_{0,1}$, see 
      Example\ref{ex:lfp}\refeq{item:pos_01}.
 
      One then uses flatness of $M$ to obtain 
      the desired vertex $\komp{b}{n}$ of a cone as follows:
      put $b_i=a_i+a'_i$ for all $i\geq 0$ and define
      $n_i:b_i\to b_{i+1}\smash b_{i+1}$ to be the one given by
      the bijection
      $$
      \Pos_{0,1}(b_i,b_{i+1}\smash b_{i+1})
      =
      \Pos_{0,1}(a_i+a'_i,b_{i+1}\smash b_{i+1})
      \cong
      \Pos_{0,1}(a_i,b_{i+1})\times\Pos_{0,1}(a'_i,b_{i+1})
      $$
      applied to the obvious pair of morphisms
      $a_i\to b_{i+1}$, $a'_i\to b_{i+1}$.
\item $\Complex{M}$ has cones for parallel pairs.

      This follows immediately from the following claim:
      \begin{itemize}
      \item[] There are no serially commutative squares
              \begin{equation}
              \label{eq:coequalized}
              \vcenter{
              \xymatrix{
              X
              \ar @<1ex> [0,1]^{u}
              \ar @<-1ex> [0,1]_-{d}
              \ar[1,0]_{s}
              &
              Y
              \ar[1,0]^{r}
              \\
              Z\smash Z
              \ar @<1ex> [0,1]^{h\smash h}
              \ar @<-1ex> [0,1]_-{l\smash l}
              &
              W\smash W
              }
              }
              \end{equation}
              whenever the maps $u$, $d$ cannot be coequalized.
      \end{itemize}
      Notice first that both $h\smash h$ and $l\smash l$
      map the ``middle element'' $(1,0)$ of $Z\smash Z$
      to the respective ``middle element'' in $W\smash W$.

      Next notice that the only reason for which $u$ and $d$ cannot
      be coequalized is that some $x\in X$ is sent to $0$ by $d$ 
      and to $1$ by $u$. Fix this $x$, and notice that equations
      $ru(x)=1$ and $rd(x)=0$ hold.

      Notice also that 
      $$
      H_h=\{ z\in Z\smash Z\mid (h\smash h)(z)=1\}
      $$ 
      is a proper subset of $\{ z\in Z\smash Z\mid z\geq m\}$
      where $m$ denotes the ``middle element'' of $Z\smash Z$.

      Similarly,
      $$
      H_l=\{ z\in Z\smash Z\mid (l\smash l)(z)=0\}
      $$ 
      is a proper subset of $\{ z\in Z\smash Z\mid z\leq m\}$.

      Especially, $H_h\cap H_l=\emptyset$.

      Suppose that the diagram~\refeq{eq:coequalized} serially
      commutes. Then $s(x)\in H_h\cap H_l$, a contradiction. 
\end{enumerate}
\end{example}

In the above example we exploited the existence of binary 
products in $\A=(\Pos_{0,1})_\fp^\op$ to observe that one can 
construct cones for two-element diagrams in $\Complex{M}$.
This is a general fact as the next result shows.

\begin{proposition}
\label{prop:nonempty=>cofiltered}
Suppose $\A$ has nonempty finite limits. 
Then the category $\Complex{M}$ is cofiltered.
\end{proposition}
\begin{proof}
Due to Assumption~\ref{ass:nonempty}, the empty
diagram in $\Complex{M}$ has a cone.

Suppose that 
$$
D:\kat{D}\to\Complex{M}
$$
with $\kat{D}$ finite and nonempty, is given. 
Let us put
$$
\xymatrix{
Dd
\ar[1,0]^{D\delta}
&
\ar @{} [1,0]|{=}
&
\dots
\ar[0,1]|-{\object @{/}}^-{m_3^d}
&
a_2^d
\ar[1,0]^{\delta_2}
\ar[0,1]|-{\object @{/}}^-{m_2^d}
&
a_1^d
\ar[1,0]^{\delta_1}
\ar[0,1]|-{\object @{/}}^-{m_1^d}
&
a_0^d
\ar[1,0]^{\delta_0}
\\
Dd'
&
&
\dots
\ar[0,1]|-{\object @{/}}_-{m_3^{d'}}
&
a_2^{d'}
\ar[0,1]|-{\object @{/}}_-{m_2^{d'}}
&
a_1^{d'}
\ar[0,1]|-{\object @{/}}_-{m_1^{d'}}
&
a_0^{d'}
\\
}
$$
and observe that, for each $n\geq 0$, its $n$-th
coordinate provides us with a diagram of shape $\kat{D}$
in $\A$. Since $\A$ has finite nonempty limits, 
we can denote, for each $n\geq 0$, by 
$$
c_n^d:a_n\to a_n^d
$$
the limit of the $n$-th coordinate.

For each $n\geq 0$, we define $m_{n+1}\in M(a_{n+1},a_n)$
as follows: since 
$$
M(a_{n+1},a_n)
\cong
\lim_d M(a_{n+1},a_n^d)
$$
holds by flatness of $M$, 
there is a unique $m_{n+1}$ such that the square
$$
\xymatrix{
a_{n+1}
\ar @{.>} [0,1]|-{\object @{/}}^-{m_{n+1}}
\ar[1,0]_{c_{n+1}^d}
&
a_n
\ar[1,0]^{c_n^d}
\\
a_{n+1}^d
\ar[0,1]|-{\object @{/}}_-{m_{n+1}^d}
&
a_n^d
}
$$
commutes.

The complex $\komp{a}{m}$ defined in the above manner is easily
seen to be a limit of $D:\kat{D}\to\Complex{M}$. 
This finishes the proof that $\Complex{M}$ is cofiltered,
hence $(\A,M)$ satisfies the Strong Sovability Condition.
\end{proof}

The Strong Solvability Condition requires, by 
Remark~\ref{rem:SSC=elementary}, the category $\Complex{M}$
to be nonempty. This is no restriction
as the following lemma shows.

\begin{lemma}
\label{lem:nonemptyness}
Either there exists no coalgebra for $M\tensor {-}$
or the category $\Complex{M}$ is nonempty.
\end{lemma}
\begin{proof}
Suppose that $e:X\to M\tensor X$ is some coalgebra. 
The functor $X$ must be flat, hence there exists an 
element $x_0\in Xa_0$.
Consider the element $e_{a_0}(x_0)\in (M\tensor X)(a_0)$. Since
$$
(M\tensor X)(a_0)=\int^a M(a,a_0)\times Xa
$$
there exist $a_1$, $m_1\in M(a_1,a_0)$ and $x_1\in Xa_1$
such that the pair $(m_1,x_2)$ represents $e_{a_0}(x_0)$. 
It is clear that in this way we can construct a complex, 
a contradiction. 
\end{proof}

\begin{definition}
\label{def:resolution}
The complex $\komp{a}{m}$ together with the sequence
$(x_n)$ constructed in the above proof is called an
{\em $e$-resolution\/} of $x_0\in Xa_0$.
\end{definition}

\begin{remark}
\label{rem:resolution}
The above construction of an $e$-resolution indicates
that a coalgebra $e:X\to M\tensor X$ is a system of 
recursive equations that ``varies in time''. For 
at ``time'' $a_0$ we can write the system of formal
recursive equations
\begin{eqnarray*}
x_0
&
\equiv
&
m_1\tensor x_1
\\
x_1
&
\equiv
&
m_2\tensor x_2
\\
&
\vdots
&
\end{eqnarray*}
where $(x_n)$ and $\komp{a}{m}$ form the $e$-resolution
of $x_0\in Xa_0$. Above, we use the tensor notation
to denote, e.g., by $m_1\tensor x_1$ the element
of $\int^a M(a,a_0)\times Xa$ represented by the
pair $(m_1,x_1)$.

Of course, any ``evolution of time''
$f:a_0\to a'_0$ provides us with a compatible corresponding
recursive system starting at $x'_0=Xf(x_0)\in Xa'_0$. 
\end{remark}

\begin{assumption}
\label{ass:nonempty}
We assume further on that $\Complex{M}$ is a nonempty category.
\end{assumption}

\begin{remark}
The proof of the following theorem is a straightforward modification
of the proof of Theorem~5.11 of~\cite{leinster1}. The reason is that 
our definition of the carrier of the final coalgebra 
(as a certain colimit) coincides with the definition of 
Tom Leinster's (as being pointwise a set of connected components
of a certain diagram, see Theorem~2.1 of~\cite{pare}). 
Observing this, the reasoning of
the proof goes exactly as in~\cite{leinster1}. 
\end{remark}

\begin{theorem}
\label{th:strong=>final}
Any $(\A,M)$ satisfying the Strong Solvability Condition
admits a final coalgebra for $M\tensor {-}$.
\end{theorem}
\begin{proof}
Define $I:\A\to\Set$ to be the colimit of the diagram
\begin{equation}
\label{eq:important_diagram}
\vcenter{
\xymatrix{
\Bigl(\Complex{M}\Bigr)^\op
\ar[0,1]^-{\pr_0^\op}
&
\A^\op
\ar[0,1]^-{Y}
&
\Presh{\A}
}
}
\end{equation}
By the Strong Solvability Condition, $I$ is a flat functor, being a
filtered colimit of representables. 
Observe that $x\in Ia$ is an equivalence class of complexes 
of the form
$$
\xymatrix{
\dots
\ar[0,1]|-{\object @{/}}^-{m_3}
&
a_2
\ar[0,1]|-{\object @{/}}^-{m_2}
&
a_1
\ar[0,1]|-{\object @{/}}^-{m_1}
&
a_0=a
}
$$
where two such complexes are equivalent if and only if there
is a zig-zag of complex morphisms having identity on $a$ as 
the $0$-th component. Thus it is exactly the description 
of elements of a final coalgebra
that Tom Leinster has for his setting in~\cite{leinster1}, page~25.
We denote equivalence classes by square brackets.

We define the coalgebra structure $\iota:I\to M\tensor I$ objectwise. 
For each $a\in \A$
$$ 
\iota_{a} :Ia\to (M\tensor I)(a) =\int^{a'} M(a',a)\times Ia' 
$$ 
is a map sending the equivalence class
$$
[
\xymatrix{
\dots
\ar[0,1]|-{\object @{/}}^-{m_3}
&
a_2
\ar[0,1]|-{\object @{/}}^-{m_2}
&
a_1
\ar[0,1]|-{\object @{/}}^-{m_1}
&
a_0=a
}
]
$$
to the element
$$
\xymatrix{
a_1
\ar[0,1]|-{\object @{/}}^-{m_1}
&
a_0
}
\tensor 
[
\xymatrix@1{
\dots
\ar[0,1]|-{\object @{/}}^-{m_3}
&
a_2
\ar[0,1]|-{\object @{/}}^-{m_2}
&
a_1
}
]
$$
of $(M\tensor I)(a)$ (recall the tensor notation of
Remark~\ref{rem:resolution}).

By Proposition~5.8 of~\cite{leinster1} such $\iota$ is 
a natural isomorphism. That $\iota:I\to M\tensor I$
is a final coalgebra follows from Theorem~5.11 of~\cite{leinster1},
once we have verified that $I$ is flat. Tom Leinster proves
finality with respect to componentwise flat functors so,
{\em a fortiori\/}, the coalgebra $\iota$ is final
with respect to coalgebras whose carriers are flat functors. 
\end{proof}

\begin{remark}
Observe that (the $a$-th component of) the mapping 
$\iota_a:Ia\to (M\tensor I)(a)$
is indeed very similar to the coalgebraic structure 
$\tau=\langle\head,\tail\rangle$
of the final coalgebra of streams of Example~\ref{ex:simple}.
\end{remark}

\begin{example}
\label{ex:final_coalgebra_on_not_lfp_part3}
(Continuation of Examples~\ref{ex:final_coalgebra_on_not_lfp_part1}
and~\ref{ex:final_coalgebra_on_not_lfp_part2}.)

We indicate how the description of the final coalgebra
for the squaring functor on the category $\Pos_{0,1}$
that we gave in Example~\ref{ex:final_coalgebra_on_not_lfp_part1}
corresponds to the description given by the proof of 
Theorem~\ref{th:strong=>final}.

We denote the module, corresponding to the squaring functor
$X\mapsto X\smash X$, by $M$. Observe that 
$$
M(a,b)
=\Pos_{0,1}(b,a\smash a)
$$
holds.

Recall that by Remark~\ref{rem:fa=flat} there is an equivalence 
$$
\Pos_{0,1}
\simeq
\Flat{(\Pos_{0,1})_\fp^\op}
$$
of categories that we will use now: the flat functor
$I:(\Pos_{0,1})_\fp^\op\to\Set$ that is the carrier
of the final coalgebra for $M\tensor {-}$ is transferred
by the above equivalence to the poset
$$
\Colim{I}{E}
$$
see Remark~\ref{rem:fa=flat}.
We define now the map
$$
\beh:
\Colim{I}{E}
\to
[0,1]
$$
where $[0,1]$ is the unit interval with the coalgebra
structure described in
Example~\ref{ex:final_coalgebra_on_not_lfp_part1}.

The mapping $\beh$ assigns to the equivalence class
$$
\Bigl[ 
[\komp{a}{m}],x\in a_0
\Bigr]
\in
\Colim{I}{E}
$$
a dyadic expansion that encodes the behaviour
of $x\in a_0$ as follows: 
we know that a complex $\komp{a}{m}$ is
a chain
$$
m_1:a_0\to a_1\smash a_1,
\quad
m_2:a_1\to a_2\smash a_2,
\quad
\dots,
\quad
m_i:a_i\to a_{i+1}\smash a_{i+1},
\quad
\dots
$$
of morphisms in $\Pos_{0,1}$. 
The morphism $m_1$ sends $x$ to the left-hand copy or 
to the right-hand copy of $a_1$, so it gives rise 
to a binary digit $k_1\in \{0,1\}$ and a new element 
$x_1\in a_1$. (If a $m_1(x)$ is in the glueing of 
the two copies of $a_1$, choose 0 or 1 arbitrarily). 
Iterating gives a binary representation $0.k_1 k_2 \dots$ 
of an element of $[0,1]$.

We will prove that $\beh$ is well-defined and a bijection. 
\begin{enumerate}
\item $\beh$ is well-defined: 
Let 
$\Bigl[[\komp{a}{m}],x\in a_0\Bigr] =\Bigl[[\komp{a'}{m'}],x'\in
  a'_0\Bigr]$, 
then there is an element $\Bigl[[\komp{c}{q}],y\in c_0\Bigr]$ of 
the colimit and a zig-zag:
\begin{center}
 $\xymatrix{
a_0\ar[r]^<(0.2){m_1}\ar[d] &  a_1\smash a_1\ar[d] \\
c_0\ar[r]^<(0.3){q_1} &  c_1\smash c_1 \\
a'_0\ar[r]^<(0.3){m'_1}\ar[u] &  a'_1\smash a'_1\ar[u]
}$
$\quad
 \xymatrix{
a_1\ar[r]^<(0.2){m_2}\ar[d] &  a_2\smash a_2\ar[d] & \\
c_1\ar[r]^<(0.3){q_2} &  c_2\smash c_2 & \dots\\
a'_1\ar[r]^<(0.3){m'_2}\ar[u] &  a'_2\smash a'_2\ar[u] &
}$
\end{center}
such that all the squares to be commutative. 

Observe that, in order to have the commutativity of the above 
squares, the morphisms $m_i, q_i, m'_i$, $i=1,2,\dots$ must 
have the same ``behaviour''. This means that if, e.g., 
the morphism $m_1$ sends $x$ to the left-hand copy of 
$a_1\smash a_1$ then also the $q_1,m'_1$ will send the 
corresponding elements to the left-hand copy of 
$c_1\smash c_1$ and $ a'_1\smash a'_1$ respectively. 
So, we take the same binary representation in $[0,1]$, i.e.,
the equality
$$
\beh([\komp{a}{m}],x\in a_0) =\beh([\komp{a'}{m'}],x'\in b_0)
$$
holds.
\item $\beh$ is one to one: \\
The key-point here is that there is a morphism 
$f: 5\to 5\smash 5$ , where 5 is the linear order with 
five elements, such that for each 
$m_i :a_i\to a_{i+1}\smash a_{i+1}$ there is a commutative square
$$ 
\xymatrix{
a_i\ar[r]^<(0.2){m_i}\ar[d]_{h} 
&  
a_{i+1}\smash a_{i+1}\ar[d]^{h'\smash h'} 
\\
5\ar[r]_<(0.3)f 
&  
5\smash 5
}
$$   
Suppose that $\{0,t_1, t_2, t_3, 1\}$ are the elements of 5, 
then the elements of $5\smash 5$ will be denoted by
$\{0,t_{1}^L, t_{2}^L, t_{3}^L, c', t_{1}^R, t_{2}^R, t_{3}^R, 1\}$. 

We define:
\begin{center}
$ 
f(t) =
\left\{ 
\scriptstyle{\begin{array}{rl}
                             0, & \mbox{if $t=0$} \\
                             1, & \mbox{if $t=1$} \\
                       t_{2}^L, & \mbox{if $t =t_1$} \\
                       t_{2}^R, & \mbox{if $t =t_3$} \\
                             c',& \mbox{if $t =t_2$}
                            \end{array}}
\right.
\quad
h(x) =\left\{ 
\scriptstyle{\begin{array}{rl}
                             0, & \mbox{if $m_{i}(x)=0$} \\
                             1, & \mbox{if $m_{i}(x)=1$} \\
                         t_{1}, & \mbox{if $m_{i}(x)\in a_{i+1}^L$} \\
                         t_{3}, & \mbox{if $m_{i}(x)\in a_{i+1}^R$} \\
                         t_{2}, & \mbox{if $m_{i}(x)=c$}  
                            \end{array}}
\right.
\quad
h'(z) =
\left\{ 
\scriptstyle{\begin{array}{rl}
                             0, & \mbox{if $z=0$} \\
                             1, & \mbox{if $z=1$} \\
                         t_{2}, & \mbox{otherwise} 
                             \end{array}}
\right.
$
\end{center}
where  $L,R$ denotes the left-hand and the right-hand copy 
and $c,c'$ are the glueing points of 
$a_{i+1}\smash a_{i+1}$ and $5\smash 5$, respectively. 
From the above it is easy to verify the commutativity of the square.

Now, if $\beh([\komp{a}{m}],x\in a_0) =\beh([\komp{b}{n}],y\in b_0)$, 
i.e., if the binary representations are the same, 
without loss of generality we can choose the 
$m_i$ and $n_i$ to send the $x_i$, $y_i$ to the same copy 
left-hand or right-hand, respectively. 
(Hence we avoid the case one of them sending an element 
to the glueing point). 
Using commutativity of the above square we have that all 
the following squares commute:
\begin{center}
 $\xymatrix{
a_0\ar[r]^<(0.2){m_1}\ar[d]_{h} &  a_1\smash a_1\ar[d]^{h'\smash h'} \\
5\ar[r]^<(0.3)f &  5\smash 5 \\
b_0\ar[r]^<(0.3){n_1}\ar[u]^{h} &  b_1\smash b_1\ar[u]_{h'\smash h'}
}$
$\quad
 \xymatrix{
a_1\ar[r]^<(0.2){m_2}\ar[d]_{h} &  a_2\smash a_2\ar[d]^{h'\smash h'} & \\
5\ar[r]^<(0.3)f &  5\smash 5 & \dots\\
b_1\ar[r]^<(0.3){n_2}\ar[u]^{h} &  b_2\smash b_2\ar[u]_{h'\smash h'} &
}$
\end{center}
From this we deduce that there is a zig-zag between the 
two complexes,$\komp{a}{m}$,$\komp{b}{n}$. 
Therefore, the equality
$$ 
[\komp{a}{m}] = [\komp{b}{n}] 
$$
holds.
\item $\beh$ is {\em epi}: 
For each binary representation $0.k_1 k_2 \dots$ of 
an element of $[0,1]$ we can find an element of the colimit, 
using the three-element linear order 3, and a sequence 
$$  
m_1:3\to 3\smash 3,
\quad
m_2:3\to 3\smash 3,
\quad
\dots
\quad
m_i:3\to 3\smash 3,
\quad
\dots
$$
of morphisms, where each $m_i$ assigns the 
middle element of 3, to the middle element in the 
left-hand copy of $3\smash 3$ if $k_i =0$, 
or the middle element in the right-hand copy if $k_i =1$.
\end{enumerate} 
\end{example}

In the realm of locally finitely presentable categories, 
{\em every\/} finitary endofunctor admits a final coalgebra.
The well-known technique for proving this result is that of
2-categorical limits of locally finitely presentable categories,
see, e.g., \cite{mp} or~\cite{ar}.

Our technique will allow us to give an alternative proof
of this theorem, see Corollary~\ref{cor:barr} below. In fact,
the colimit of~\refeq{eq:important_diagram} gives an explicit 
description of a final coalgebra.

\begin{corollary}
\label{cor:barr}
Every finitary endofunctor of a locally finitely presentable category
admits a final coalgebra.
\end{corollary}
\begin{proof}
Recall that the category of the form $\Flat{\A}$ 
is locally finitely presentable, if the category
$\A$ has {\em all\/} finite limits.
Denote by $(\A,M)$ the corresponding self-similarity system. 
We need to show that $\Complex{M}$ is cofiltered.
\begin{enumerate}
\item 
The category $\Flat{\A}\simeq\Lex{\A}$ has an initial object, $\bot$,
say. Hence the unique morphism $!:\bot\to M\tensor \bot$ is a
coalgebra and the category $\Complex{M}$
is nonempty by Lemma~\ref{lem:nonemptyness}.
\item
By Proposition~\ref{prop:nonempty=>cofiltered}, the category
$\Complex{M}$ has cones for nonempty finite diagrams. 
\end{enumerate}
Now use Theorem~\ref{th:strong=>final}.
\end{proof}

Theorem~\ref{th:strong=>final} provides us with a concrete
description of the final coalgebra as the colimit
of the filtered diagram 
$$
\xymatrix{
\Bigl(\Complex{M}\Bigr)^\op
\ar[0,1]^-{\pr_0^\op}
&
\A^\op
\ar[0,1]^-{Y}
&
\Presh{\A}
}
$$ 
From that one can easily deduce, for example, the well-known description
of the final coalgebra for the endofunctor $X\mapsto X\times X+A$
on $\Set$ that we gave in the Introduction.

\section{The Weak Solvability Condition}
\label{sec:WSC}

Cofilteredness of the category $\Complex{M}$ may be hard
to verify in the absence of finite limits in $\A$.
We give here a weaker condition that is easier to verify.
In particular, we are going to replace the Strong Solvability 
Condition by a condition of the same type
but ``holding just on the head of complexes''.
This whole section is devoted to finding conditions
of ``how to propagate from the head of a complex 
to the whole complex''.
Proving the existence of a final coalgebra 
will require though some extra finiteness condition 
on the module $M$, see Definition~\ref{def:ess_finite}.
Our condition is a weakening of that considered by Tom
Leinster~\cite{leinster1} in connection with self-similar
objects in topology. 
The main result of this section, Theorem~\ref{th:weak=>strong}, 
then shows that this finiteness condition allows us to conclude
that a final coalgebra exists.
Our argument applies to self-similarity systems considered
by Tom Leinster~\cite{leinster1} and therefore strenghtens his
result on the existence of final coalgebras for self-similarity
systems.

The key tool for the propagation technique
is ``K\"{o}nig's Lemma for preorders'', see
Theorem~\ref{th:koenig_for_posets} below.
The result relies on a topological fact proved
by Arthur Stone in~\cite{stone}.

To be able to state the weak condition we first need to generalize 
filteredness of a category to filteredness of a functor.

\begin{definition}
\label{def:filtered_functor}
A functor $F:\kat{X}\to\kat{Y}$ is called {\em filtering\/},
if there exists a cocone for the composite $F\o D$,
for every functor $D:\kat{D}\to\kat{X}$ with $\kat{D}$ finite.

A functor $F$ is called {\em cofiltering\/} if $F^\op$ is
filtering.
\end{definition}

\begin{remark}
Hence a category $\kat{X}$ is filtered if and only if
the identity functor $\Id:\kat{X}\to\kat{X}$ is filtering.
\end{remark}

A natural candidate for a weaker form of solvability
condition is the following one.

\begin{definition}
\label{def:WSC}
We say that $(\A,M)$ satisfies the
{\em Weak Solvability Condition\/} if the functor 
$$
\pr_0:\Complex{M}\to\nComplex{M}{0}
$$ 
is cofiltering.
\end{definition}

In particular, observe that the Weak Solvability Condition
holds when the category $\A$ is cofiltered.

\begin{remark}
\label{rem:WSC=elementary}
In elementary terms, the Weak Solvability Condition 
says the following three conditions:
\begin{enumerate}
\item The category $\A$ is non-empty.
\item For every pair $\komp{a}{m}$,
      $\komp{a'}{m'}$ in $\Complex{M}$ there is a span
      $$
      \xymatrixrowsep{.5pc}
      \xymatrix{
      &
      a_0
      \\
      b
      \ar[-1,1]^{f}
      \ar[1,1]_{f'}
      &
      \\
      &
      a'_0
      }
      $$
      in $\A$.
\item For every parallel pair of the form
      $$
      \xymatrix{
      \komp{a}{m}
      \ar @<1ex> [0,1]^-{(u_\bullet)}
      \ar @<-1ex> [0,1]_-{(v_\bullet)}
      &
      \komp{a'}{m'}
      }
      $$
      in $\Complex{M}$
      there is a fork
      $$
      \xymatrix{
      b
      \ar[0,1]^-{f}
      &
      a_0
      \ar @<1ex> [0,1]^-{u_0}
      \ar @<-1ex> [0,1]_-{v_0}
      &
      a'_0
      }
      $$
      in $\A$.
\end{enumerate}
Observe that, since we assume that $\Complex{M}$ is nonempty
(Assumption~\ref{ass:nonempty}), the above condition~(1) is 
satisfied: the category $\A$ is nonempty.
\end{remark}

Observe that if $(\A,M)$ satisfies the Strong Solvability Condition, 
it does satisfy the Weak Solvability Condition. In fact, in this
case {\em every\/} functor $\pr_n:\Complex{M}\to\nComplex{M}{n}$
is cofiltering. The following result shows that the Weak Solvability
Condition can be formulated in this way.

\begin{proposition}
\label{prop:WSC}
The following are equivalent:
\begin{enumerate}
\item The Weak Solvability Condition.
\item The functors $\pr_n:\Complex{M}\to\nComplex{M}{n}$
      are cofiltering for all $n\geq 0$.
\end{enumerate}
\end{proposition}
\begin{proof}
That~(2) implies~(1) is clear. To prove the converse, we need
to verify the following three properties:
\begin{enumerate}
\item[(a)] Every category $\nComplex{M}{n}$ is non-empty. This is clear: 
      we assume that that $\Complex{M}$ is non-empty, see Assumption~\ref{ass:nonempty}.
\item[(b)] Every pair $\pr_n\komp{a}{m}$, $\pr_n\komp{a'}{m'}$ in
      $\nComplex{M}{n}$ has a cone.

      Observe that, due to Weak Solvability Condition applied at stage
      $n$, we have the following diagram
      $$
      \xymatrixrowsep{.5pc}
      \xymatrix{
      & 
      a_n\ar[0,1]|-{\object @{/}}^-{m_n} 
      & 
      a_{n-1}\ar[0,1]|-{\object @{/}}^-{m_{n-1}} 
      & \dots\ar[0,1]|-{\object @{/}}^-{m_1} 
      & a_0 
      \\
      b_n
      \ar@{.>}[ur]^{f_n}
      \ar@{.>}[dr]_{f'_n} 
      & 
      & 
      & 
      & 
      \\
      & 
      a'_n\ar[0,1]|-{\object @{/}}_-{m'_n} 
      & 
      a'_{n-1}\ar[0,1]|-{\object @{/}}_-{m'_{n-1}} 
      & 
      \dots\ar[0,1]|-{\object @{/}}_-{m'_1} 
      & 
      a'_0   
      }
      $$
      Since the functor $M(b_n,{-})$ is flat, the pair
      $m_n @ f_n\in M(b_n,a_{n-1})$, 
      $m'_n @ f'_n\in M(b_n,a'_{n-1})$
      of its elements has a cone: 
      $$
      \xymatrixrowsep{.5pc}
      \xymatrix{
      & 
      a_n\ar[0,1]|-{\object @{/}}^-{m_n} 
      & a_{n-1} 
      \\
      b_n
      \ar@{.>}[ur]^{f_n}
      \ar@{.>}[dr]_{f'_n}
      \ar@{.>}[0,1]|-{\object @{/}}^-{u_n} 
      & 
      b_{n-1}
      \ar@{.>}[ur]
      \ar@{.>}[dr] 
      &  
      \\
      & 
      a'_n
      \ar[0,1]|-{\object @{/}}_-{m'_n} 
      & 
      a'_{n-1}
      }
      $$
      If we proceed like this down to zero we obtain the desired
      vertex $\komp{b}{u}^{(n)}$ in $\nComplex{M}{n}$:
      $$
      \xymatrixrowsep{.5pc}
      \xymatrix{
      &
      \pr_n\komp{a}{m} 
      \\
      \komp{b}{u}^{(n)}
      \ar[-1,1]^{(f_\bullet)}
      \ar[1,1]_{(f'_\bullet)} 
      & 
      \\
      & 
      \pr_n\komp{a'}{m'}
      }
      $$
\item[(c)] For every parallel pair of the form
      $$
      \xymatrix{
      \pr_n\komp{a}{m}
      \ar @<1ex> [0,1]^-{\pr_n(u_\bullet)}
      \ar @<-1ex> [0,1]_-{\pr_n(v_\bullet)}
      &
      \pr_n\komp{a'}{m'}
      }
      $$
      in $\nComplex{M}{n}$, there is a fork.
      
      Consider the following diagram:
      $$
      \xymatrix{
      b_{n}
      \ar@{.>}[d]_{f_n}
      \ar@{.>}[0,1]|-{\object @{/}}^{l_n} 
      & 
      b_{n-1}
      \ar@{.>}[d]_{f_{n-1}}
      \ar@{.>}[0,1]|-{\object @{/}}^{l_{n-1}} 
      &
      \dots
      \ar@{.>}[0,1]|-{\object @{/}}^{l_1}  
      & 
      b_0\ar@{.>}[d]_{f_0}
      \\
      a_n
      \ar @<-1ex> [d]_-{u_n}
      \ar @<1ex> [d]^-{v_n}
      \ar[0,1]|-{\object @{/}}^{m_n} 
      & 
      a_{n-1}
      \ar @<-1ex> [d]_-{u_{n-1}}
      \ar @<1ex> [d]^-{v_{n-1}}
      \ar[0,1]|-{\object @{/}}^{m_{n-1}} 
      &
      \dots\ar[0,1]|-{\object @{/}}^{m_1} 
      & 
      a_0
      \ar @<-1ex> [d]_-{u_0}
      \ar @<1ex> [d]^-{v_0} 
      \\
      a'_n 
      \ar[0,1]|-{\object @{/}}_{m'_n} 
      &
      a'_{n-1}\ar[0,1]|-{\object @{/}}_{m'_{n-1}} 
      &
      \dots\ar[0,1]|-{\object @{/}}_{m'_1} 
      & 
      a_0
      }
      $$
      Again, start at stage $n$, use the Weak Solvability Condition
      there to obtain $f_n$, and then use flatness of $M(b_n,{-})$
      to obtain $l_n$ and $f_{n-1}$. Proceed like this down to zero
      and obtain the desired fork
      $$
      \xymatrix{
      \komp{b}{l}^{(n)}
      \ar@{.>}[0,1]^-{(f_\bullet)^{(n)}} 
      & 
      \pr_n\komp{a}{m}
      \ar @<1ex> [0,1]^-{\pr_n(u_\bullet)}
      \ar @<-1ex> [0,1]_-{\pr_n(v_\bullet)}
      & 
      \pr_n\komp{a'}{m'}
      }
      $$
      in $\nComplex{M}{n}$.
\end{enumerate}
This finishes the proof.
\end{proof}

In the proof that the Weak Solvability Condition implies
the Strong one, we will need to use ``K\"{o}nig's Lemma'' for 
preorders that we formulate in Theorem~\ref{th:koenig_for_posets}
below. 

Recall that a {\em preorder\/}
$\strukt{X,\less}$ is a set $X$ equipped with a reflexive,
transitive binary relation $\less$. 

Recall also that a subset $B\subseteq X$ of a preorder 
is called {\em downward-closed\/},
if for every $b\in B$ and $b'\less b$ we have $b'\in B$.
The dual notion is called {\em upward-closed\/}.

A subset $S$ of a preorder $\strukt{X,\less}$ is called
{\em final\/} if for every $x\in X$ there exists $s\in S$
with $x\less s$.

\begin{theorem}
\label{th:koenig_for_posets}
Suppose that 
\begin{equation}
\label{eq:chain}
\vcenter{
\xymatrix{
\dots
\ar[0,1]
&
\kat{P}_{n+1}
\ar[0,1]^-{p^{n+1}_{n}}
&
\kat{P}_n
\ar[0,1]^-{p^{n}_{n-1}}
&
\dots
\ar[0,1]^-{p^{1}_{0}}
&
\kat{P}_0
}
}
\end{equation}
is a chain of preorders and monotone maps, 
that satisfies the following two conditions:
\begin{enumerate}
\item Every $\kat{P}_n$ has a nonempty finite final subset.
\item The image of any upward-closed set under 
      $p^{n+1}_{n}:\kat{P}_{n+1}\to\kat{P}_n$
      is upward-closed.
\end{enumerate}
Then the limit $\lim\kat{P}_n$ is nonempty, i.e., there
is a sequence $(x_n)$ with $p^{n+1}_{n}(x_{n+1})=x_n$
holding for every $n\geq 0$.
\end{theorem}

The proof of Theorem~\ref{th:koenig_for_posets} will rely
on some facts from General Topology that we recall now.
As a reference to topology we refer to the book~\cite{engelking}.

Recall that every preorder $\strukt{X,\less}$ can be 
equipped with the {\em lower topology\/} $\tau_\less$, 
if we declare the open sets to be exactly the downward closed 
sets.

Observe that a set $B$ is closed in the topology $\tau_\less$
if and only if it is upward-closed.

\medskip\noindent
{\it Proof of Theorem~\ref{th:koenig_for_posets}.}
The assumptions~(1) and~(2) of the statement of the
theorem assure that each $\kat{P}_n$ is a nonempty compact
space in its lower topology and each $p^{n+1}_n$
is a closed continuous map (i.e., on top of continuity, the image of 
a closed set is a closed set). By result of Arthur
Stone~\cite{stone}, Theorem~2, 
any $\omega^\op$-chain of nonempty compact
spaces and closed continuous maps has a nonempty
limit. Therefore $\lim\kat{P}_n$ is nonempty.
\mbox{}\hfill\qed

\begin{remark}
Of course, Theorem~\ref{th:koenig_for_posets} holds whenever 
Conditions~(1) and~(2)
hold ``cofinally'', i.e., whenever there exists $n_0$ such that
Conditions~(1) and~(2) hold for all $n\geq n_0$.
\end{remark}

\begin{notation}
\label{not:chain}
For any diagram $D:\kat{D}\to\Complex{M}$ with $\kat{D}$ finite, 
let $\kat{P}^D_n$ denote the following preorder:
\begin{enumerate}
\item Points of $\kat{P}^D_n$ are cones for the composite
      $\pr_n\o D:\kat{D}\to\nComplex{M}{n}$.
\item The relation $c\less_n c'$ holds in $\kat{P}^D_n$ 
      if and only if the cone $c$ factors through the cone $c'$.
\end{enumerate}
For each $n\geq 0$ denote by
$$
p^{n+1}_n:\kat{P}^D_{n+1}\to\kat{P}^D_n
$$
the obvious restriction map and observe that it is monotone.
\end{notation}

Also observe that the Weak Solvability Condition guarantees that every
preorder $\kat{P}^D_n$ is nonempty by Proposition~\ref{prop:WSC}.
The Weak Solvability Condition alone does not imply 
the Strong one --- the self-similarity system $(\A,M)$ has to fulfill 
additional conditions that will allow us to apply
Theorem~\ref{th:koenig_for_posets}. 

\begin{definition}
\label{def:ess_finite}
We say that the module $M$ is {\em compact\/}, 
if the preorder $\kat{P}^D_n$ has a nonempty finite final subset, 
for each $n\geq 0$ and each finite nonempty diagram 
$D:\kat{D}\to\Complex{M}$.
\end{definition}

First we give easy examples of compact modules.

\begin{example}
\mbox{}\hfill
\begin{enumerate}
\item 
Every module $M$ on a finitely complete category $\A$ is compact:
in fact, in this case every preorder $\kat{P}_n^D$ has a one-element
final set.
\item
If the module $M$ is finite in the sense of~\cite{leinster1}, i.e.,
if every functor $M({-},b):\A^\op\to\Set$ has a finite category
of elements, then it is compact.
\end{enumerate}
\end{example}

Nontrivial examples of compact modules will follow later from
Proposition~\ref{prop:em=>compact}, see Example~\ref{ex:linear}. 
We need to recall the 
concept of a factorization system for cocones first. For details,
see, e.g., Chapter~IV of~\cite{ahs}.

\begin{definition}
Let $\kat{K}$ be a finitely accessible category.
\begin{enumerate}
\item We say that a cocone $c_d:Dd\to X$ is {\em jointly epi\/} if,
      for every parallel pair $u$, $v$,
      the equality $u\o c_d=v\o c_d$ for all $d$
      implies that $u=v$ holds.
\item We say that $\kat{K}$ is a 
      {\em (finite jointly epi, extremal mono)-category\/}
      if the following two conditions are satisfied:
      \begin{enumerate}
      \item Every cocone $c_d:Dd\to X$ for a finite diagram
            can be factored as
            $$
            \xymatrix{
            Dd
            \ar[0,1]^-{e_d}
            &
            Z
            \ar[0,1]^-{j}
            &
            X
            }
            $$
            where $e_d$ is jointly epi and $j$ is extremal
            mono.
      \item For every commutative square
            $$
            \vcenter{
            \xymatrix{
            Dd
            \ar[0,1]^-{e_d}
            \ar[1,0]_{f_d}
            &
            X
            \ar[1,0]^{g}
            \\
            A
            \ar[0,1]_-{j}
            &
            B
            }
            }
            \quad
            \mbox{(for all $d$)}
            $$
            where $e_d$ is a jointly epi cocone and $j$ is extremal
            mono, there is a unique diagonal $m:X\to A$ making the
            obvious triangles commutative.
      \end{enumerate}
\item We say that $\kat{K}_\fp$ is {\em finitely cowellpowered\/},
      if every finite diagram $D:\kat{D}\to\kat{K}_\fp$ admits
      (up to isomorphism) only a nonempty finite set of jointly 
      epi cocones.
\end{enumerate}
\end{definition}

\begin{proposition}
\label{prop:em=>compact}
Suppose the finitely accessible category $\kat{K}$ satisfies
the following conditions:
\begin{enumerate}
\item  $\kat{K}$ is a (finite jointly epi, extremal mono)-category.
\item $\kat{K}_\fp$ is finitely cowellpowered.
\end{enumerate}
Suppose that a finitary functor $\Phi:\kat{K}\to\kat{K}$
preserves extremal monos. Then the flat module corresponding to
$\Phi$ is compact.
\end{proposition}
\begin{proof}
We will use the description of complexes from 
Remark~\ref{rem:complexes=simple}. 

Let $D:\kat{D}\to\Complex{M}$ be a finite nonempty
diagram. Choose any $n\geq 0$ and
denote the value of the composite
$\pr_n\o D$ by commutative squares
$$
\xymatrix{
\pr_n\o Dd
\ar[1,0]^-{\pr_n\o D\delta}
&
\ar @{} [1,0]|{=}
&
a_0^d
\ar[0,1]^-{m_1^d}
&
\Phi (a_1^d)
&
a_1^d
\ar[0,1]^-{m_2^d}
&
\Phi (a_2^d)
\ar @{} [1,2]|{\objectstyle\dots}
&
&
a_{n-1}^d
\ar[0,1]^-{m_n^d}
&
\Phi (a_n^d)
\\
\pr_n\o Dd'
&
&
a_0^{d'}
\ar[0,1]_-{m_0^{d'}}
\ar[-1,0]^{\delta_0}
&
\Phi (a_1^{d'})
\ar[-1,0]_{\Phi(\delta_1)}
&
a_1^{d'}
\ar[0,1]_-{m_2^{d'}}
\ar[-1,0]^{\delta_1}
&
\Phi (a_2^{d'})
\ar[-1,0]_{\Phi(\delta_2)}
&
&
a_{n-1}^{d'}
\ar[0,1]_-{m_n^{d'}}
\ar[-1,0]^{\delta_{n-1}}
&
\Phi (a_n^{d'})
\ar[-1,0]_{\Phi(\delta_n)}
}
$$ 
in $\kat{K}$. We will construct the finite nonempty
{\em initial\/} (notice the change of the 
variance: $\A$ is $\kat{K}_\fp^\op$)
of cocones for $\pr_n\o D$ by proceeding from
$i=n-1$ downwards to $0$ 
as follows:
\begin{enumerate}
\item[] For every jointly epi
      cocone $e_{i+1}:a_{i+1}^d\to z_{i+1}$
      choose all jointly epi cocones $e_i^d:a_i^d\to z_i$ and 
      all connecting morphisms $c_{i+1}:z_i\to\Phi(z_{i+1})$ 
      making the following diagram
      $$
      \xymatrix{
      z_i
      \ar @{.>} [0,1]^-{c_{i+1}}
      &
      \Phi(z_{i+1})
      \\
      a_i^d
      \ar @{.>} [-1,0]^{e_i^d}
      \ar[0,1]_-{m_{i+1}^d}
      &
      \Phi(a_{i+1}^d)
      \ar[-1,0]_{\Phi(e_{i+1}^d)}
      }
      $$
      commutative. Observe that there is at least one 
      such pair: the factorization of the cocone 
      $\Phi(e_{i+1}^d)\o  m_{i+1}^d$ into a jointly
      epi and extremal mono. Since every cocone 
      $e_i^d$ is jointly epi, the corresponding $c_{i+1}$
      is determined uniquely. 
\end{enumerate}
We claim that the above nonempty finite family of cocones
for $\pr_n\o D$ is initial. To that end, consider any 
cocone
$$
\xymatrix{
w_0
\ar[0,1]^-{f_1}
&
\Phi (w_1)
&
w_1
\ar[0,1]^-{f_2}
&
\Phi (w_2)
\ar @{} [1,2]|{\objectstyle\dots}
&
&
w_{n-1}
\ar[0,1]^-{f_n}
&
\Phi (w_n)
\\
a_0^{d}
\ar[0,1]_-{m_0^{d}}
\ar[-1,0]^{g^d_0}
&
\Phi (a_1^{d})
\ar[-1,0]_{\Phi(g^d_1)}
&
a_1^{d}
\ar[0,1]_-{m_2^{d}}
\ar[-1,0]^{g^d_1}
&
\Phi (a_2^{d})
\ar[-1,0]_{\Phi(g^d_2)}
&
&
a_{n-1}^{d}
\ar[0,1]_-{m_n^{d}}
\ar[-1,0]^{g^d_{n-1}}
&
\Phi (a_n^{d})
\ar[-1,0]_{\Phi(g^d_n)}
}
$$  
for $\pr_n\o D$. Factorize the cocone $g^d_n$
into a jointly epi $e^d_n:a^d_n\to z_n$
followed by an extremal mono $j_n:z_n\to w_n$.
Do the same thing for the cocone 
$g^d_{n-1}$ and then use
the diagonalization property to obtain the
desired $c_{n}:z_{n-1}\to \Phi(z_n)$
$$
\xymatrix{
w_{n-1}
\ar[0,1]^-{f_n}
&
w_n
\\
z_{n-1}
\ar @{.>} [0,1]^-{c_n}
\ar[-1,0]^{j_{n-1}}
&
\Phi(z_{n})
\ar[-1,0]_{\Phi(j_n)}
\\
a_{n-1}^d
\ar[-1,0]^{e_{n-1}^d}
\ar[0,1]_-{m_{n}^d}
&
\Phi(a_{n}^d)
\ar[-1,0]_{\Phi(e_{n}^d)}
}
$$
using the fact that $\Phi(j_n)$ is extremal mono
by assumption. Proceed like this downwards to $0$
and obtain thus one of the above chosen cocones
through which the given cocone of $g$ factorizes. 
\end{proof}

\begin{example}
\label{ex:linear}
Recall from Example~\ref{ex:lfp}\refeq{item:linear} that the category
$\Lin$ of all linear orders and all monotone maps is finitely
accessible. We indicate that it fulfills the assumptions
of the above proposition and give several examples of finitary
endofunctors that preserve extremal monos.
\begin{enumerate}
\item Jointly epi cocones $e_d:Dd\to X$ are exactly thoses
      where (the underlying set of) $X$ is the union of the
      images of all $Dd$.
\item A monotone map $j:A\to B$ is an extremal mono
      if and only if $j$ is injective and the linear order
      on $A$ is that induced by $B$.
\end{enumerate}
From the above it is clear that $\Lin$ is 
a (finite jointly epi, extremal mono)-category and that
$\Lin_\fp$ is finitely cowellpowered.

To give various examples of functors that preserve 
extremal monos, we need to introduce the following 
notation: given linear orders $X$ and $Y$ we denote by
$$
X\after Y
\quad
\mbox{(read: $X$ {\em then\/} $Y$)}
$$ 
the linear order on the disjoint union of (the underlying sets of) 
$X$ and $Y$ by putting every element of $X$ to be lower than any
element of $Y$ and leaving the linear orders of $X$ and $Y$ unchanged.

The second construction is that of {\em ordinal product\/}, by
$$
X*Y
$$ 
we denote the linear order on the cartesian product of 
(underlying sets of) $X$ and $Y$ where we replace
each element of $Y$ by a disjoint copy of $X$.
More precisely, $(x,y)<(x',y')$ holds if and only if
either $x<x'$ holds or $x=x'$ and $y<y'$.

It can be proved easily that, for example, the following 
two assignments
$$
X\mapsto X*\omega,
\quad
X\mapsto (X*\omega)\after 1
$$
where $\omega$ is the first countable ordinal and $1$ denotes
the one-element linear order, are finitary functors and they
both preserve extremal monos.
\end{example}

Our main result on compact modules is the following one.

\begin{theorem}
\label{th:weak=>strong}
Suppose that $M$ is a compact module. Then
the Weak Solvability Condition implies the Strong one. 
\end{theorem}
\begin{proof}
We know that $\Complex{M}$ is nonempty by Assumption~\ref{ass:nonempty}.
We have to construct a cone for every diagram
$D:\kat{D}\to\Complex{M}$ with $\kat{D}$
finite nonempty. 

Form the corresponding chain
\begin{equation}
\label{eq:chain2}
\vcenter{
\xymatrix{
\dots
\ar[0,1]
&
\kat{P}^D_{n+1}
\ar[0,1]^-{p^{n+1}_{n}}
&
\kat{P}^D_n
\ar[0,1]^-{p^{n}_{n-1}}
&
\dots
\ar[0,1]^-{p^{1}_{0}}
&
\kat{P}^D_0
}
}
\end{equation}
of preorders and monotone maps. We will verify first that it satisfies
Conditions~(1) and~(2) of Theorem~\ref{th:koenig_for_posets}.
\begin{enumerate}
\item Each $\kat{P}^D_n$ contains a nonempty finite final subset
      since the module $M$ is assumed to be compact.
\item The image of every upward-closed set under the monotone map 
      $p^{n+1}_{n}$ is upward-closed.

      Denote the value of $D:\kat{D}\to\Complex{M}$ by
      $$
      \xymatrix{
      Dd
      \ar[1,0]^{D\delta}
      &
      \ar @{} [1,0]|{=}
      &
      \dots
      \ar[0,1]|-{\object @{/}}^-{m_3^d}
      &
      a_2^d
      \ar[1,0]^{\delta_2}
      \ar[0,1]|-{\object @{/}}^-{m_2^d}
      &
      a_1^d
      \ar[1,0]^{\delta_1}
      \ar[0,1]|-{\object @{/}}^-{m_1^d}
      &
      a_0^d
      \ar[1,0]^{\delta_0}
      \\
      Dd'
      &
      &
      \dots
      \ar[0,1]|-{\object @{/}}_-{m_3^{d'}}
      &
      a_2^{d'}
      \ar[0,1]|-{\object @{/}}_-{m_2^{d'}}
      &
      a_1^{d'}
      \ar[0,1]|-{\object @{/}}_-{m_1^{d'}}
      &
      a_0^{d'}
      \\
      }
      $$
      Then the value of $\pr_n\o D:\kat{D}\to\nComplex{M}{n}$ is given
      by
      $$
      \xymatrix{
      \pr_n\o Dd
      \ar[1,0]^{\pr_n\o D\delta}
      &
      \ar @{} [1,0]|{=}
      &
      a_n^d
      \ar[0,1]|-{\object @{/}}^-{m_n^d}
      \ar[1,0]^{\delta_n}
      &
      \dots
      \ar[0,1]|-{\object @{/}}^-{m_3^d}
      &
      a_2^d
      \ar[1,0]^{\delta_2}
      \ar[0,1]|-{\object @{/}}^-{m_2^d}
      &
      a_1^d
      \ar[1,0]^{\delta_1}
      \ar[0,1]|-{\object @{/}}^-{m_1^d}
      &
      a_0^d
      \ar[1,0]^{\delta_0}
      \\
      \pr_n\o Dd'
      &
      &
      a_n^{d'}
      \ar[0,1]|-{\object @{/}}_-{m_n^{d'}}
      &
      \dots
      \ar[0,1]|-{\object @{/}}_-{m_3^{d'}}
      &
      a_2^{d'}
      \ar[0,1]|-{\object @{/}}_-{m_2^{d'}}
      &
      a_1^{d'}
      \ar[0,1]|-{\object @{/}}_-{m_1^{d'}}
      &
      a_0^{d'}
      \\
      }
      $$ 
      for every $n\geq 0$.
      
      Choose an upward-closed set $S\subseteq \kat{P}^D_{n+1}$. 
      Every $s\in S$ is a cone
      for the above diagram $\pr_{n+1}\o D$ and we denote this cone by
      $$
      s=
      \vcenter{
      \xymatrix{
      s_{n+1}
      \ar[1,0]^{\sigma^d_{n+1}}
      \ar[0,1]|-{\object @{/}}^-{m_{n+1}^s}
      &
      \dots
      \ar[0,1]|-{\object @{/}}^-{m_{3}^s}
      &
      s_{2}
      \ar[1,0]^{\sigma^d_{2}}
      \ar[0,1]|-{\object @{/}}^-{m_{2}^s}
      &
      s_{1}
      \ar[1,0]^{\sigma^d_{1}}
      \ar[0,1]|-{\object @{/}}^-{m_{1}^s}
      &
      s_{0}
      \ar[1,0]^{\sigma^d_{0}}
      \\
      a_{n+1}^d
      \ar[0,1]|-{\object @{/}}_-{m_{n+1}^d}
      &
      \dots
      \ar[0,1]|-{\object @{/}}_-{m_3^d}
      &
      a_2^d
      \ar[0,1]|-{\object @{/}}_-{m_2^d}
      &
      a_1^d
      \ar[0,1]|-{\object @{/}}_-{m_1^d}
      &
      a_0^d
      }
      }
      $$
      Choose any $s$ in $S$ and consider $b$ in $\kat{P}^D_n$
      such that $p^{n+1}_{n}(s)\less_n b$ holds.
      We need to find $s\less_{n+1} t$ such that
      $p^{n+1}_n(t)=b$.     

      In our notation, $b$ has the form
      $$
      b=
      \vcenter{
      \xymatrix{
      b_{n}
      \ar[1,0]^{\beta^d_{n}}
      \ar[0,1]|-{\object @{/}}^-{m_{n}^b}
      &
      \dots
      \ar[0,1]|-{\object @{/}}^-{m_{3}^b}
      &
      b_{2}
      \ar[1,0]^{\beta^d_{2}}
      \ar[0,1]|-{\object @{/}}^-{m_{2}^b}
      &
      b_{1}
      \ar[1,0]^{\beta^d_{1}}
      \ar[0,1]|-{\object @{/}}^-{m_{1}^b}
      &
      b_{0}
      \ar[1,0]^{\beta^d_{0}}
      \\
      a_{n}^d
      \ar[0,1]|-{\object @{/}}_-{m_{n}^d}
      &
      \dots
      \ar[0,1]|-{\object @{/}}_-{m_3^d}
      &
      a_2^d
      \ar[0,1]|-{\object @{/}}_-{m_2^d}
      &
      a_1^d
      \ar[0,1]|-{\object @{/}}_-{m_1^d}
      &
      a_0^d
      }
      }
      $$

      The inequality $p^{n+1}_{n}(s)\less_n b$ means that 
      there exists a diagram of the form
      $$
      \xymatrix{
      s_{n}
      \ar[1,0]^{g_{n}}
      \ar[0,1]|-{\object @{/}}^-{m_{n}^s}
      &
      \dots
      \ar[0,1]|-{\object @{/}}^-{m_{3}^s}
      &
      s_{2}
      \ar[1,0]^{g_{2}}
      \ar[0,1]|-{\object @{/}}^-{m_{2}^s}
      &
      s_{1}
      \ar[1,0]^{g_{1}}
      \ar[0,1]|-{\object @{/}}^-{m_{1}^s}
      &
      s_{0}
      \ar[1,0]^{g_{0}}
      \\
      b_{n}
      \ar[1,0]^{\beta^d_{n}}
      \ar[0,1]|-{\object @{/}}^-{m_{n}^b}
      &
      \dots
      \ar[0,1]|-{\object @{/}}^-{m_{3}^b}
      &
      b_{2}
      \ar[1,0]^{\beta^d_{2}}
      \ar[0,1]|-{\object @{/}}^-{m_{2}^b}
      &
      b_{1}
      \ar[1,0]^{\beta^d_{1}}
      \ar[0,1]|-{\object @{/}}^-{m_{1}^b}
      &
      b_{0}
      \ar[1,0]^{\beta^d_{0}}
      \\
      a_{n}^d
      \ar[0,1]|-{\object @{/}}_-{m_{n}^d}
      &
      \dots
      \ar[0,1]|-{\object @{/}}_-{m_3^d}
      &
      a_2^d
      \ar[0,1]|-{\object @{/}}_-{m_2^d}
      &
      a_1^d
      \ar[0,1]|-{\object @{/}}_-{m_1^d}
      &
      a_0^d
      }
      $$      
      where the equalities
      $\beta^d_i\o g_i=\sigma^d_i$ hold for every
      $i\in\{ 0,\dots,n\}$. 

      Consider the following diagram:
      $$
      \xymatrix{
      s_{n+1}
      \ar@{=} [1,0]
      \ar[0,1]|-{\object @{/}}^-{m_{n+1}^s}
      &
      s_{n}
      \ar[1,0]^{g_{n}}
      \ar[0,1]|-{\object @{/}}^-{m_{n}^s}
      &
      \dots
      \ar[0,1]|-{\object @{/}}^-{m_{3}^s}
      &
      s_{2}
      \ar[1,0]^{g_{2}}
      \ar[0,1]|-{\object @{/}}^-{m_{2}^s}
      &
      s_{1}
      \ar[1,0]^{g_{1}}
      \ar[0,1]|-{\object @{/}}^-{m_{1}^s}
      &
      s_{0}
      \ar[1,0]^{g_{0}}
      \\
      s_{n+1}
      \ar[1,0]^{\beta^d_{n+1}}
      \ar[0,1]|-{\object @{/}}^-{g_n\o m_{n+1}^s}
      &
      b_{n}
      \ar[1,0]^{\beta^d_{n}}
      \ar[0,1]|-{\object @{/}}^-{m_{n}^b}
      &
      \dots
      \ar[0,1]|-{\object @{/}}^-{m_{3}^b}
      &
      b_{2}
      \ar[1,0]^{\beta^d_{2}}
      \ar[0,1]|-{\object @{/}}^-{m_{2}^b}
      &
      b_{1}
      \ar[1,0]^{\beta^d_{1}}
      \ar[0,1]|-{\object @{/}}^-{m_{1}^b}
      &
      b_{0}
      \ar[1,0]^{\beta^d_{0}}
      \\
      a_{n+1}^d
      \ar[0,1]|-{\object @{/}}_-{m_{n+1}^d}
      &
      a_{n}^d
      \ar[0,1]|-{\object @{/}}_-{m_{n}^d}
      &
      \dots
      \ar[0,1]|-{\object @{/}}_-{m_3^d}
      &
      a_2^d
      \ar[0,1]|-{\object @{/}}_-{m_2^d}
      &
      a_1^d
      \ar[0,1]|-{\object @{/}}_-{m_1^d}
      &
      a_0^d
      }
      $$ 
      Thus, the desired $t$ has the form
      $$
      \xymatrix{
      s_{n+1}
      \ar[1,0]^{\beta^d_{n+1}}
      \ar[0,1]|-{\object @{/}}^-{g_n @ m_{n+1}^s}
      &
      b_{n}
      \ar[1,0]^{\beta^d_{n}}
      \ar[0,1]|-{\object @{/}}^-{m_{n}^b}
      &
      \dots
      \ar[0,1]|-{\object @{/}}^-{m_{3}^b}
      &
      b_{2}
      \ar[1,0]^{\beta^d_{2}}
      \ar[0,1]|-{\object @{/}}^-{m_{2}^b}
      &
      b_{1}
      \ar[1,0]^{\beta^d_{1}}
      \ar[0,1]|-{\object @{/}}^-{m_{1}^b}
      &
      b_{0}
      \ar[1,0]^{\beta^d_{0}}
      \\
      a_{n+1}^d
      \ar[0,1]|-{\object @{/}}_-{m_{n+1}^d}
      &
      a_{n}^d
      \ar[0,1]|-{\object @{/}}_-{m_{n}^d}
      &
      \dots
      \ar[0,1]|-{\object @{/}}_-{m_3^d}
      &
      a_2^d
      \ar[0,1]|-{\object @{/}}_-{m_2^d}
      &
      a_1^d
      \ar[0,1]|-{\object @{/}}_-{m_1^d}
      &
      a_0^d
      }
      $$ 
      Hence the image of every upward-closed set under
      the monotone map $p^{n+1}_{n}:\kat{P}^D_{n+1}\to\kat{P}^D_n$
      is upward-closed.
\end{enumerate}
Therefore, by Theorem~\ref{th:koenig_for_posets},
we have an element $(x_n)$ of the limit $\lim\kat{P}^D_n$.

Denote every $x_n$ as follows:
$$
x_n=
\vcenter{
\xymatrix{
x^{n}_{n}
\ar[1,0]^{\chi^{n,d}_{n}}
\ar[0,1]|-{\object @{/}}^-{m_{n}^{x_n}}
&
\dots
\ar[0,1]|-{\object @{/}}^-{m_{3}^{x_n}}
&
x^{n}_{2}
\ar[1,0]^{\chi^{n,d}_{2}}
\ar[0,1]|-{\object @{/}}^-{m_{2}^{x_n}}
&
x^{n}_{1}
\ar[1,0]^{\chi^{n,d}_{1}}
\ar[0,1]|-{\object @{/}}^-{m_{1}^{x_n}}
&
x^{n}_{0}
\ar[1,0]^{\chi^{n,d}_{0}}
\\
a_{n}^d
\ar[0,1]|-{\object @{/}}_-{m_{n}^d}
&
\dots
\ar[0,1]|-{\object @{/}}_-{m_3^d}
&
a_2^d
\ar[0,1]|-{\object @{/}}_-{m_2^d}
&
a_1^d
\ar[0,1]|-{\object @{/}}_-{m_1^d}
&
a_0^d
}
}
$$
From that we can define a complex
$$
\xymatrixcolsep{3pc}
\xymatrix{
\dots
\ar[0,1]|-{\object @{/}}^-{m_4^{x_4}}
&
x^{3}_{3}
\ar[0,1]|-{\object @{/}}^-{m_3^{x_3}}
&
x^{2}_{2}
\ar[0,1]|-{\object @{/}}^-{m_2^{x_2}}
&
x^{1}_{1}
\ar[0,1]|-{\object @{/}}^-{m_1^{x_1}}
&
x^{0}_{0}
}
$$
that is obviously a vertex of a cone for $D:\kat{D}\to\Complex{M}$.
\end{proof}

\begin{corollary}
Every compact module satisfying the Weak Solvability
Condition has a final coalgebra.
\end{corollary}

\begin{example}
\label{ex:pavlovic}
Recall from Example~\ref{ex:linear} that the modules
corresponding to the finitary endofunctors
$$
X\mapsto X*\omega,
\quad
X\mapsto (X*\omega)\after 1
$$
of the category $\Lin$ are compact. Since $\Lin$
satisfies the Weak Solvability Conditions, the above
two functors have final coalgebras by the above corollary.
The linear orders of these coalgebras are the continuum
and Cantor space, respectively, see~\cite{pavlovic+pratt}
for a proof. 
\end{example}

\section{What the Existence of a Final Coalgebra Entails}
\label{sec:final=>WSC}

We show in this section that the existence of final coalgebras
entails the Weak Solvability Condition, provided the module
is pointed. As a corollary, we derive a necessary condition
on the category $\A$ so that the identity functor on $\Flat{\A}$
admits a final coalgebra, see Corollary~\ref{cor:identity}.

\begin{assumption}
\label{ass:pointed}
We assume that in this section that $M$ is {\em pointed\/},
i.e., that $M$ is equipped with a module morphism 
$c:\A\to M$.
\end{assumption}

Of course, the assumption is clearly satisfied if and
only if, when passing from $M$ to the finitary endofunctor $\Phi$, 
there exists a natural transformation $\Id\to\Phi$. 

\begin{remark}
From the Assumption~\ref{ass:pointed} it follows that 
every representable functor $\A(a,{-})$ admits
a coalgebra structure 
$$ 
c_a:\A(a,{-})\to M(a,{-})
$$
for $M\tensor {-}$ (we used that 
$\Bigl(M\tensor\A\Bigr)(a,{-})\cong M(a,{-})$ holds). 
This of course entails that $\Complex{M}$ is 
nonempty, see Lemma~\ref{lem:nonemptyness}.

Moreover, for every $f:a\to a'$, the natural transformation
$\A(f,{-}):\A(a',{-})\to\A(a,{-})$ is a coalgebra morphism, i.e.,
the square
\begin{equation}
\label{eq:naturality}
\vcenter{
\xymatrix{
\A(a',{-})
\ar[0,2]^-{c_{a'}}
\ar[1,0]_{\A(f,{-})}
&
&
M(a',{-})
\ar[1,0]^{M(f,{-})}
\\
\A(a,{-})
\ar[0,2]_-{c_a}
&
&
M(a,{-})
}
}
\end{equation}
commutes. 
\end{remark}

\begin{theorem}
\label{th:final=>weak}
Suppose that $M$ is pointed and suppose that a final coalgebra for
$M\tensor {-}$ exists. Then $\pr_0$ is cofiltering,
i.e., the Weak Solvability Condition holds.
\end{theorem}
\begin{proof}
Let us denote by $j:J\to M\tensor J$ the final coalgebra for
$M\tensor {-}$.

Denote by $\sol{c_a}:\A(a,{-})\to J$
the unique coalgebra morphism such that the square
$$
\xymatrix{
\A(a,{-})
\ar[0,2]^-{c_a}
\ar[1,0]_{\sol{c_a}}
&
&
M\tensor\A(a,{-})
\ar[1,0]^{M\circ\sol{c_a}}
\\
J
\ar[0,2]_-{j}
&
&
M\tensor J
}
$$
commutes.

Then the following triangle
$$
\xymatrixrowsep{1pc}
\xymatrix{
\A(a',{-})
\ar[1,1]^{\sol{c_{a'}}}
\ar[2,0]_{\A(f,{-})}
&
\\
&
J
\\
\A(a,{-})
\ar[-1,1]_{\sol{c_a}}
}
$$
commutes by finality of $j:J\to M\tensor J$ and the 
square~\refeq{eq:naturality}.

Recall that, in any case, one can form a colimit $I$
of the diagram
$$
\xymatrix{
\Bigl(\Complex{M}\Bigr)^\op
\ar[0,1]^-{\pr_0^\op}
&
\A^\op
\ar[0,1]^-{Y}
&
\Presh{\A}
}
$$
We do not claim that $I:\A\to\Set$ is flat. In fact, we
will just use the fact that $I$ is a colimit. For observe
that so far we have proved that the collection of
morphisms 
$$
\sol{c_{\pr_0^\op\komp{a}{m}}}:\A(a_0,{-})\to J
$$
forms a cocone for the diagram $Y\o\pr_0^\op$.
Hence there exists a natural transformation 
$$
\overline{\beta}:I\to J
$$ 
The natural transformation $\overline{\beta}$ induces a functor
$F:\Complex{M}\to\elts{J}$ by putting
$$
\komp{a}{m}
\mapsto
x\in Ja_0
$$
where the element $x\in Ja_0$ corresponds to the natural
transformation 
$\sol{c_{\pr_0^\op\komp{a}{m}}}:\A(a_0,{-})\to J$
by Yoneda Lemma.

Then the diagram
$$
\xymatrix{
\Complex{M}
\ar[0,2]^-{F}
\ar[1,1]_{\pr_0}
&
&
\elts{J}
\ar[1,-1]^{{\mathsf{proj}}}
\\
&
\A
&
}
$$
commutes. Since $J$ is a flat functor, the category $\elts{J}$ 
is cofiltered. Hence $\pr_0={\mathsf{proj}}\o F$ is a cofiltering 
functor.
\end{proof}

\begin{corollary}
\label{cor:identity}
If the identity functor on the category $\Flat{\A}$ has a final
coalgebra, then the category $\A$ must be cofiltered.
\end{corollary}

\begin{remark}
The above Corollary shows that the identity endofunctor
of a Scott complete category $\kat{K}$, see 
Example~\ref{ex:lfp}\refeq{item:pos_01}, canot have a final
coalgebra unless the category $\kat{K}$ is in fact
locally finitely presentable.
\end{remark}

What we have proved so far, allows us to go in full circle:

\begin{corollary}
Suppose that 
$
M:\A
\xymatrixcolsep{1.5pc}\xymatrix{\ar[0,1]|-{\object @{/}}&}
\A
$
is a pointed, compact module.
Then the following are equivalent:
\begin{enumerate}
\item The self-similarity system $(\A,M)$ satisfies the 
      Weak Solvability Condition.
\item The self-similarity system $(\A,M)$ satisfies the
      Strong Solvability Condition.
\item The colimit of the diagram
      $$
      \xymatrix{
      \Bigl(\Complex{M}\Bigr)^\op
      \ar[0,1]^-{\pr_0^\op}
      &
      \A^\op
      \ar[0,1]^-{Y}
      &
      \Presh{\A}
      }
      $$
      is a flat functor.
\item The final coalgebra for $M\tensor {-}$ exists.
\end{enumerate}
\end{corollary}

\section{Conclusions and Future Research}
\label{sec:conclusions}

We have provided a new uniform way of constructing
final coalgebras for finitary endofunctors of locally
finitely presentable categories. We have argued about
the necessity of expanding these results to the case
of finitely accessible categories. To that end we have
formulated general conditions that are sufficient for
the existence of a final coalgebra. We expect that our
conditions can be exploited for finding new interesting
examples of final coalgebras in accessible categories.

In many concrete examples where the final coalgebra cannot
exist for cardinality reasons (e.g., the categories where
all maps are injections) we expect that suitable modifications
of our results will provide coalgebras of ``rational terms''.
This means coalgebras comprising of solutions of finitary
recursive systems, see~\cite{mccoy}.


\begin{thebibliography}{MMMM}
\bibitem[AAMV]{aamv}
        P.~Aczel, J.~Ad\'{a}mek, S.~Milius and J.~Velebil,
        Infinite trees and completely iterative theories:
        A coalgebraic view, {\em Theoret.~Comput.~Sci.} 300 (2003),
        1--45
\bibitem[A${}_1$]{adamek}
        J.~Ad\'{a}mek,
        On final coalgebras of continuous functors,
        {\em Theoret.~Comput.~Sci.} 294 (2003), 3--29
\bibitem[A${}_2$]{adamek_scc}
        J.~Ad\'{a}mek,
        A categorical generalization of Scott domains,
        {\em Math. Structures Comput. Sci.\/} 7 (1997), 419--443
\bibitem[ABLR]{ablr}
        J.~Ad\'amek, F.~Borceux, S.~Lack and J.~Rosick\'y,
        A classification of accessible categories,
        {\em Jour. Pure Appl. Alg.\/} 175 (2002), 7--30
\bibitem[AHS]{ahs}
        J.~Ad\'amek, H.~Herrlich and G.~Strecker,
        {\em Abstract and concrete categories\/},
        John Wiley \& Sons, New York, 1990,
        available electronically as a TAC reprint No.~17 at
        {\tt http://www.tac.mta.ca/tac/reprints/articles/17/tr17abs.html}
\bibitem[AMV]{mccoy}
        J.~Ad\'{a}mek, S.~Milius and J.~Velebil,
        Iterative algebras at work,
        {\em Math. Structures Comput. Sci.\/} 16 (2006), 1085--1131
\bibitem[ADJ]{adj}
        J.~A.~Goguen, S.~W.~Thatcher, E.~G.~Wagner and J.~B.~Wright,
        Initial algebra semantics and continuous algebras,
        {\em Journal ACM} 24 (1977), 68--95
\bibitem[ARu]{america+rutten}
        P.~America and J.~J.~M.~M.~Rutten, Solving reflexive
        domain equations in a category of complete
        metric spaces, {\em J. Comput. System Sci.\/}
        39 (1989), 343--375
\bibitem[AR]{ar}
        J.~Ad\'amek and J.~Rosick\'y,
        {\em Locally presentable and accessible categories},
        Cambridge University Press, 1994  
\bibitem[Bo]{borceux}
        F.~Borceux,
        {\em Handbook of categorical algebra\/}
        (three volumes), Cambridge University Press, 1994
\bibitem[D${}_1$]{diers1}
        Y.~Diers, 
        Cat\'{e}gories multialg\'{e}briques,
        {\em Arch. Math. (Basel)\/} 34 (1980), 193--209
\bibitem[D${}_2$]{diers2}
        Y.~Diers, 
        Cat\'{e}gories localement multipr\'{e}sentables,
        {\em Arch. Math. (Basel)\/} 34 (1980), 344--356
\bibitem[E]{elgot}
        C.~C.~Elgot, Monadic computation and iterative algebraic
        theories, in: {\em Logic Colloquium `73\/}
        (eds: H.~E.~Rose and J.~C.~Shepherdson), North-Holland
        Publishers, Amsterdam, 1975
\bibitem[EBT]{ebt}
        C.~C.~Elgot, S.~L.~Bloom and R.~Tindell,
        On the algebraic structure of rooted trees,
        {\em J. Comput. System Sci.\/} 16 (1978), 361--399
\bibitem[En]{engelking}
        R.~Engelking, {\em General topology\/},
        Sigma Series in Pure Mathematics, Berlin Heldermann, 1989
\bibitem[F]{freyd}
        P.~Freyd, 
        {\em Real coalgebra\/}, posting to category theory mailing list,
        22 Dec 1999, available electronically at
        {\tt http://www.mta.ca/~cat-dist/catlist/1999/realcoalg}
\bibitem[GU]{gu}
        P.~Gabriel and F.~Ulmer,
        {\em Lokal pr\"{a}sentierbare Kategorien\/}, Lecture Notes
        in Mathematics 221, Springer 1971
\bibitem[L]{lair}
        C.~Lair,
        Cat\'{e}gories modelables et cat\'{e}gories esquissables,
        {\em Diagrammes\/} 6 (1981), L1--L20
\bibitem[Le${}_1$]{leinster1}
        T.~Leinster,
        A general theory of self-similarity I,
        {\em arXiv:math/0411344v1\/}
\bibitem[Le${}_2$]{leinster2}
        T.~Leinster, 
        A general theory of self-similarity II: recognition,
        {\em arXiv:math/0411345v1\/}
\bibitem[McL]{maclane}
        S.~MacLane,
        {\em Categories for the working mathematician\/},
        Springer Verlag, 1971
\bibitem[MPa]{mp}
        M.~Makkai and R.~Par\'{e},
        {\em Accessible categories: The foundations of categorical
          model theory\/}, Contemporary Mathematics~104, American
        Mathematical Society, 1989
\bibitem[PP]{pavlovic+pratt}
        D.~Pavlovi\'{c} and V.~Pratt,
        The continuum as a final coalgebra,
        {\em Theoret. Comput. Sci.\/} 280 (2002), 105--122
\bibitem[Pa]{pare}
        R.~Par\'{e},
        Connected components and colimits,
        {\em Jour. Pure Appl. Alg.\/} 3 (1973), 21--42
\bibitem[R]{rutten}
        J.~J.~M.~M.~Rutten, Universal coalgebra: A theory of
        systems, {\em Theoret. Comput. Sci.\/} 249 (2000), 3--80
\bibitem[S]{stone}
        A.~H.~Stone,
        Inverse limits of compact spaces.
        {\em General Topology Appl.\/} 10 (1979), no. 2, 203--211.
\end{thebibliography}
\end{document}